\documentclass[11pt,a4paper]{article}
\usepackage[T1]{fontenc}
\usepackage[utf8]{inputenc}
\usepackage{lmodern}
\usepackage{amsmath,amssymb,amsthm,mathtools,bm}
\usepackage[margin=28mm,top=25mm,bottom=27mm]{geometry}
\usepackage{microtype}
\usepackage[hidelinks]{hyperref}
\numberwithin{equation}{section}
\newtheorem{theorem}{Theorem}[section]
\newtheorem{proposition}[theorem]{Proposition}
\newtheorem{lemma}[theorem]{Lemma}
\theoremstyle{remark}

\allowdisplaybreaks[2]
\hypersetup{pdftitle={A Proof of a Polya-Szego Conjecture: A PDE Approach to Capacity Inequalities},pdfauthor={Zikang Deng}}
\title{A Proof of a P\'olya--Szeg\H{o} Conjecture\\[0.3em]
\large A PDE Approach to Capacity Inequalities}
\author{Zikang Deng\\\normalsize Beijing Normal University}
\date{October 1, 2026}
\begin{document}
\maketitle

\begin{abstract}
We prove a planar capacity-ratio inequality conjectured by P\'olya and Szeg\H{o}, together with its extension to adjacent dimensions: for every integer $n\geq2$ and compact set $K\subset\mathbb R^n$,
\[
C_{n-1}(K)\leq
\left[\frac{\Gamma(n/2)}{\sqrt\pi\,\Gamma((n+1)/2)}\right]^{1/(n-1)}C_{n-2}(K),
\]
where $C_s$ denotes the length-homogeneous Riesz capacity and, when $n=2$, $C_{n-2}=C_0$ is logarithmic capacity. Closed balls attain equality, so the constant is sharp. Starting from an exterior Dirichlet problem and its nonnegative boundary flux, we construct a function $\rho$ in one higher dimension, with coordinates $X=(x,z)$, and prove that, for $z\ne0$,
\[
\Delta\rho=\frac{n-1}{\rho}|\nabla\rho|^2+\frac{\rho_z}{z},
\qquad
\rho|\nabla\rho|^2\leq(\rho^2-1)\frac{\rho_z}{z}.
\]
A decreasing function $\Phi_n$ solving a suitable ordinary differential equation makes $U=\Phi_n(\rho)$ superharmonic outside the conductor. Green's formula on the exterior of level sets of $U$ then bounds the Dirichlet integral by the total flux at infinity, and the Dirichlet principle yields the capacity bound. Smooth outer approximation extends the result to arbitrary compact sets. The same construction, combined with a comparison principle for weighted elliptic equations, gives the sharp upper bound for $C_q/C_{n-2}$ when $n-2<q<n$.
\end{abstract}
\noindent\textbf{Keywords:} Riesz capacity; logarithmic capacity; exterior Dirichlet problem; superharmonic function; gradient estimate; Dirichlet principle.

\section{Introduction}
Let $n\geq2$, let $K\subset\mathbb R^n$ be a nonempty compact set, and denote by $\mathcal P(K)$ the set of Borel probability measures supported on $K$. For $0<s<n$, define the Riesz energy and capacity by
\begin{equation}\label{eq:1.1}
\begin{aligned}
I_s(\mu)&=\iint |x-y|^{-s}\,d\mu(x)\,d\mu(y),\\
E_s(K)&=\inf_{\mu\in\mathcal P(K)}I_s(\mu),
&C_s(K)&=E_s(K)^{-1/s}.
\end{aligned}
\end{equation}
In the plane, we also define logarithmic energy and logarithmic capacity by
\begin{equation}\label{eq:1.2}
\begin{aligned}
I_0(\mu)&=\iint\log\frac1{|x-y|}\,d\mu(x)\,d\mu(y),\\
E_0(K)&=\inf_{\mu\in\mathcal P(K)}I_0(\mu),
&C_0(K)&=e^{-E_0(K)}.
\end{aligned}
\end{equation}
The logarithmic kernel is bounded below on the Cartesian square of a fixed compact set, so $E_0(K)$ cannot be $-\infty$. When the energy is $+\infty$, the corresponding capacity is defined to be zero; the capacity of the empty set is also defined to be zero. This normalization gives capacity the dimension of length. More precisely, for $r>0$ and $\mu\in\mathcal P(K)$, let $\mu_r$ be the pushforward of $\mu$ under $x\mapsto rx$, so that $\mu_r(A)=\mu(r^{-1}A)$. Substitution in the kernels gives
\[
I_s(\mu_r)=r^{-s}I_s(\mu)\quad(s>0),
\qquad I_0(\mu_r)=I_0(\mu)-\log r.
\]
Since $\mu\mapsto\mu_r$ is a bijection from $\mathcal P(K)$ onto $\mathcal P(rK)$, taking infima and applying the definitions yields $C_s(rK)=rC_s(K)$. Translation preserves distances, so capacity is also translation invariant. For monotonicity with respect to sets, if $K\subset L$, then $\mathcal P(K)\subset\mathcal P(L)$, whence $E_s(K)\geq E_s(L)$ and $C_s(K)\leq C_s(L)$.

For $n\geq3$, $C_{n-2}$ is the length-homogeneous form of Newtonian capacity in $\mathbb R^n$. When $K$ is embedded isometrically in the horizontal hyperplane $\mathbb R^n\times\{0\}\subset\mathbb R^{n+1}$, distances remain unchanged, while the Newtonian kernel in the higher-dimensional space becomes $|X-Y|^{-(n-1)}$, corresponding to $C_{n-1}$. We study the ratio of the capacities of the same set in these two adjacent dimensions, using logarithmic capacity in the lower dimension when $n=2$.

In 1945, P\'olya and Szeg\H{o} asked whether, among planar compact sets with prescribed logarithmic capacity, the disk has the greatest three-dimensional electrostatic capacity \cite[Conjecture (1.3)]{ref1}. In the normalization above, this asks for the inequality
\[
C_1(K)\leq\frac2\pi C_0(K),\qquad K\subset\mathbb R^2.
\]
Clark and Laugesen formulated the problem as an extremal problem for the ratio of Newtonian capacities in adjacent dimensions \cite[Conjecture 1.1]{ref2}. We prove the following results. Throughout, $\Gamma$ denotes the Gamma function,
\[
\Gamma(t)=\int_0^\infty r^{t-1}e^{-r}\,dr,\qquad t>0,
\]
and $S^{d-1}$ denotes the unit sphere in $\mathbb R^d$, with surface area $|S^{d-1}|$.

\begin{theorem}[Capacity ratio in adjacent dimensions]\label{thm:1.1}
Let $n\geq2$ be an integer. For every compact set $K\subset\mathbb R^n$,
\begin{equation}\label{eq:1.3}
C_{n-1}(K)\leq\kappa_n C_{n-2}(K),\qquad
\kappa_n=\left[\frac{\Gamma(n/2)}{\sqrt\pi\,\Gamma((n+1)/2)}\right]^{1/(n-1)}.
\end{equation}
Every closed ball of positive radius attains equality, so $\kappa_n$ is sharp.
\end{theorem}

\begin{theorem}[Capacity inequality for general exponents]\label{thm:1.2}
Let $n\geq2$ and $n-2<q<n$. For every compact set $K\subset\mathbb R^n$,
\begin{equation}\label{eq:1.4}
C_q(K)\leq A_{n,q}^{1/q}C_{n-2}(K),\qquad
A_{n,q}=\frac{\Gamma(n/2)}{\Gamma(1+q/2)\Gamma((n-q)/2)}.
\end{equation}
Every closed ball of positive radius attains equality, so the constant $A_{n,q}^{1/q}$ is sharp.
\end{theorem}

Neither result requires $K$ to be connected, convex, or smooth. Theorem~1.1 is the case $q=n-1$ of Theorem~1.2. In this case the ordinary Laplace equation and Dirichlet integral apply directly, and we therefore give its proof first.

As for previous progress, Clark and Laugesen established a conditional result in \cite[Conjecture 1.2 and Theorem 1.3]{ref2}: their proposed comparison for the vertical Dirichlet energy in a suitable strip with Neumann boundary conditions would imply the adjacent-dimensional capacity-ratio conjecture. Starting from the exterior capacitary potential, we prove a pointwise gradient estimate and use it to construct elliptic supersolutions, yielding Theorems~1.1 and~1.2. Capacity ratios with both exponents variable are also studied in \cite{ref3}. Here the smaller exponent is always fixed at $n-2$, and the larger exponent is restricted to $q\in(n-2,n)$.

\subsection{Preliminary results}
We collect the classical potential-theoretic results used in the proof and specify their normalizations. The bound on the equilibrium potential, the Dirichlet principle, and the capacity formula for closed balls are established results; their sources are given after each statement. These results will be used to construct the extension and prove the capacity inequalities. The notation needed in each result is given in its statement, and all other notation is defined at first use.

\begin{theorem}[Existence of equilibrium measures and an upper bound for equilibrium potentials]\label{thm:1.3}
For the energy in \eqref{eq:1.1} or \eqref{eq:1.2}, the energy infimum on a nonempty compact set is attained by a probability measure. A probability measure attaining this infimum is called an equilibrium measure for the corresponding energy. If $K\subset\mathbb R^n$, $n-2<q<n$, and $C_q(K)>0$, its $q$-equilibrium measure $\sigma_q$ satisfies
\begin{equation}\label{eq:1.5}
\int_K|x-y|^{-q}\,d\sigma_q(y)\leq E_q(K),\qquad x\in\mathbb R^n.
\end{equation}
\end{theorem}

Existence follows from weak compactness of probability measures and lower semicontinuity of the energy. If $E_s(K)<+\infty$, take measures with $I_s(\mu_j)\to E_s(K)$. Weak compactness gives a weakly convergent subsequence with limit $\mu\in\mathcal P(K)$. Here weak convergence means that $\int f\,d\mu_j\to\int f\,d\mu$ for every continuous function $f\in C(K)$. Since the kernel is lower semicontinuous and bounded below on $K\times K$,
\[
E_s(K)\leq I_s(\mu)\leq\liminf_j I_s(\mu_j)=E_s(K),
\]
so $\mu$ attains the infimum. The case in which every measure has infinite energy is immediate. See \cite[Lemma 4.1.3]{ref6} for this existence result and its potential-theoretic background. Inequality~\eqref{eq:1.5} is the maximum principle for Riesz equilibrium potentials in the range $0<n-q<2$; see \cite[Theorem 2.1]{ref8}, whose parameter is $\alpha=n-q\in(0,2)$. This range includes the case $n=2$ considered here. In Section~6, this whole-space bound will allow us to compare the normalized $q$-equilibrium potential with the supersolution we construct.

\begin{theorem}[Dirichlet principle and capacity normalization]\label{thm:1.4}
Let $d\geq3$ and let $E\subset\mathbb R^d$ be compact. Define
\begin{equation}\label{eq:1.6}
\operatorname{cap}_{D,d}(E)=\inf\left\{\int_{\mathbb R^d}|\nabla v|^2\,dX:
\begin{array}{l}v\in C_c^\infty(\mathbb R^d),\\ v\geq1\text{ in a neighborhood of }E\end{array}\right\}.
\end{equation}
Here $C_c^\infty$ denotes the smooth functions with compact support. Then
\begin{equation}\label{eq:1.7}
\operatorname{cap}_{D,d}(E)=(d-2)|S^{d-1}|C_{d-2}(E)^{d-2},
\end{equation}
where $C_{d-2}(E)$ is defined by \eqref{eq:1.1} in $\mathbb R^d$.
\end{theorem}

This equivalence is the Dirichlet principle for Newtonian capacity \cite{ref5}. The identity relating the Dirichlet integral of a potential to the energy of its measure, together with the normalization constants, is also given in \cite[Lemma 4.2]{ref2}.

To explain the constant in \eqref{eq:1.7}, write $\gamma_d=(d-2)|S^{d-1}|$. The fundamental solution of $-\Delta$ is $\gamma_d^{-1}|X|^{2-d}$. For a smooth conductor $E$, its capacitary potential $h$ is harmonic in the exterior, equals $1$ on the conductor, and tends to zero at infinity. Let $\nu$ be the unit normal pointing from the conductor into its exterior, and extend $h$ by $1$ inside the conductor. Green's formula gives
\[
\int_{\mathbb R^d}|\nabla h|^2\,dX
=-\int_{\partial E}\partial_\nu h\,dS
=\gamma_d C_{d-2}(E)^{d-2}.
\]
The Dirichlet principle states that the energy of $h$ is the infimum in \eqref{eq:1.6}. Although $h$ itself is not compactly supported, this value is obtained by truncating and smoothing $h$ and approximating its boundary values. The details will be given in Section~2 for the exterior capacitary potential $w$. For arbitrary compact sets, smooth outer approximation and continuity of both capacities under decreasing limits give \eqref{eq:1.7}. Continuity under decreasing limits for the energy definition is proved in Section~5. For \eqref{eq:1.6}, the same property follows directly from the infimum definition, since each admissible function is at least $1$ in a neighborhood of $E$.

\begin{theorem}[Capacity of a closed ball]\label{thm:1.5}
Let $B^n=\{x\in\mathbb R^n:|x|\leq1\}$ be the closed unit ball. Then $C_{n-2}(B^n)=1$, and, for $n-2<q<n$,
\begin{equation}\label{eq:1.8}
C_q(B^n)^q=\frac{\Gamma(n/2)}{\Gamma(1+q/2)\Gamma((n-q)/2)}.
\end{equation}
\end{theorem}

This formula is given in \cite[Appendix A]{ref4}. For $n\geq3$, the identity $C_{n-2}(B^n)=1$ also follows directly from the exterior capacitary potential $|x|^{2-n}$; for $n=2$, the corresponding exterior Green function is $\log|x|$. The sharpness of our inequalities will follow from this closed-ball formula.

\subsection{Outline of the proof}
First let $F$ be a smooth conductor with $C^\infty$ boundary and connected complement, and write $c=C_{n-2}(F)>0$. Starting from the exterior capacitary potential in $\mathbb R^n\setminus F$, we use the nonnegative measure given by its boundary normal derivative to construct a function $\rho$ on $\mathbb R^{n+1}$. We identify $F$ with $F\times\{0\}$ and write $X=(x,z)$. The resulting function equals $1$ on $F$, is strictly greater than $1$ outside $F$, and satisfies
\[
\rho(X)=\frac{|X|}{c}+O(1),\qquad |X|\to\infty.
\]
In Section~2, differential identities for the integral kernels show that $\rho$ satisfies the elliptic equation
\[
\Delta\rho=\frac{n-1}{\rho}|\nabla\rho|^2+\frac{\rho_z}{z}.
\]
Here $\Delta$ and $\nabla$ act on all the variables $X$, and $\rho_z$ is the partial derivative with respect to $z$. At points with $z=0$ and $x\notin F$, the quotient $\rho_z/z$ is interpreted by continuity.

Section~3 establishes the key pointwise gradient estimate
\[
\rho|\nabla\rho|^2\leq(\rho^2-1)\frac{\rho_z}{z}.
\]
Its proof uses the maximum principle for the exterior capacitary potential to control an integral at an arbitrary horizontal point, followed by H\"older's inequality and completion of the square. In two dimensions, Jensen's inequality for the logarithm plays the corresponding role.

We then choose a decreasing function with $\Phi_n(1)=1$ and $\Phi_n(\infty)=0$ satisfying
\[
\Phi_n''(r)+\left(\frac{n-1}{r}+\frac{r}{r^2-1}\right)\Phi_n'(r)=0.
\]
This choice allows us to control the Laplacian of $U=\Phi_n(\rho)$ by the gradient estimate, obtaining $-\Delta U\geq0$. Section~4 computes its asymptotic behavior at infinity:
\[
U(X)=a_nc^{n-1}|X|^{1-n}+O(|X|^{-n}),\qquad
a_n=\frac{\Gamma(n/2)}{\sqrt\pi\,\Gamma((n+1)/2)}.
\]
To turn the supersolution into a capacity bound, Section~5 applies Green's formula on the exterior of regular level sets of $U$, and truncates $U$ to obtain admissible functions equal to $1$ in a neighborhood of $F$. The total flux at infinity bounds their Dirichlet integrals and gives
\[
\operatorname{cap}_{D,n+1}(F)\leq(n-1)|S^n|a_nc^{n-1}.
\]
Together with \eqref{eq:1.7}, this proves the main inequality in the smooth case. Smooth outer approximation and continuity of capacity under decreasing limits then give the result for arbitrary compact sets. Finally, Section~6 replaces the Laplacian by the weighted elliptic operator associated with the Riesz kernel and applies a comparison principle to prove Theorem~1.2.

\section{Exterior capacitary potentials and extension to a higher dimension}
\label{sec:exterior}

Let $D\subset\mathbb R^n$ be a nonempty bounded open set whose boundary is a compact $C^\infty$ hypersurface, let $F=\overline D$, and suppose that the exterior domain $\Omega=\mathbb R^n\setminus F$ is connected. We allow $D$ to have finitely many connected components. We do not assume that $\Omega$ is simply connected. In what follows, $B_R$ denotes the open ball of radius $R$ centered at the origin in the Euclidean space under consideration. Let $\nu$ be the unit normal pointing from $F$ into $\Omega$. Thus the outward normal to the truncated exterior domain $B_R\setminus F$ on its inner boundary $\partial F$ is $-\nu$. Set $c=C_{n-2}(F)$, where $C_{n-2}=C_0$ when $n=2$. Since $D$ contains a closed ball of positive radius, monotonicity of capacity and Theorem~1.5 give $c>0$.

Embed $F$ in the hyperplane $z=0$ of $\mathbb R^{n+1}$, retaining the notation $F$. Write $X=(x,z)$ and $R_y=|X-(y,0)|$. Unless otherwise stated, $\nabla$ and $\Delta$ act on all $n+1$ variables, whereas $\nabla_x$ and $\Delta_x$ act only on $x$.

\subsection{The exterior Dirichlet problem in higher dimensions}

Suppose that $n\geq3$, and set $p=n-2$ and $\gamma_n=p|S^{n-1}|$. The existence and uniqueness theorem for the smooth exterior Dirichlet problem yields a function $w\in C^\infty(\Omega)\cap C(\overline\Omega)$ satisfying
\begin{equation}
\begin{cases}
\Delta_x w=0&\text{in }\Omega,\\
w=1&\text{on }\partial F,\\
w(x)\longrightarrow0&\text{as }|x|\longrightarrow\infty.
\end{cases}
\label{eq:2.1}
\end{equation}
Here the relevant exterior problem is that for a bounded $C^\infty$ conductor, which need not be connected. Boundary regularity implies that $w$ is smooth up to $\partial F$ from the exterior. Applying the maximum principle first on $B_R\setminus F$ and then letting $R\to\infty$ gives $0\leq w\leq1$. Since $\Omega$ is connected and $w$ is nonconstant, the strong maximum principle gives $0<w<1$ in $\Omega$. Extend $w$ to be $1$ in $F$; the resulting function is still continuous. The Hopf boundary point lemma gives $-\partial_\nu w>0$ on $\partial F$.

We now give the flux representation and its normalization in detail. Define the finite nonnegative measure
\begin{equation}
 d\eta(y)=\frac{-\partial_\nu w(y)}{\gamma_n}\,dS(y).
\label{eq:2.2}
\end{equation}
Its density is smooth. Denote the action of a distribution $T$ on a test function $\varphi$ by $\langle T,\varphi\rangle$. For every $\varphi\in C_c^\infty(\mathbb R^n)$, integration by parts gives
\[
\langle-\Delta_x w,\varphi\rangle
=\int_\Omega\nabla_x w\cdot\nabla_x\varphi\,dx
=-\int_{\partial F}\varphi\partial_\nu w\,dS
=\gamma_n\int_{\partial F}\varphi\,d\eta.
\]
Here $\nabla_xw=0$ in the interior of the conductor, and the normal on the inner boundary of the exterior domain is $-\nu$, which accounts for the sign of the boundary term.

Set
\[
v(x)=\int_{\partial F}|x-y|^{-p}\,d\eta(y).
\]
This integral is finite everywhere and continuous. Indeed, the local surface area of the smooth compact hypersurface satisfies
$S(\partial F\cap B(x,t))\leq Ct^{n-1}$, and the density of $\eta$ is bounded. For $0<t<1$, decompose $B(x,t)\setminus\{x\}$ into the annuli
$2^{-k-1}t<|x-y|\leq2^{-k}t$, $k=0,1,\ldots$. Then
\[
\begin{aligned}
\int_{\partial F\cap B(x,t)}|x-y|^{-p}\,d\eta(y)
&\leq C\sum_{k=0}^\infty(2^{k+1}/t)^p(2^{-k}t)^{n-1}\\
&\leq Ct\sum_{k=0}^\infty2^{-k}\leq C't.
\end{aligned}
\]
To prove continuity at a fixed point $x_0$, choose $0<t<1/3$ and $|x-x_0|<t$. The fixed near region $B(x_0,2t)$ is contained in $B(x,3t)$, so both integrals over this near region are $O(t)$ by the preceding estimate. Outside this region, the distance to the singularity of the kernel is at least $t$, and the dominated convergence theorem applies. Letting first $x\to x_0$ and then $t\downarrow0$ proves continuity of $v$. Denote by $\delta_y$ the unit point mass at $y$. The fundamental solution identity
$-\Delta_x|x-y|^{-p}=\gamma_n\delta_y$ gives
$-\Delta_xv=\gamma_n\eta=-\Delta_xw$. Thus $w-v$ is distributionally harmonic on the whole space, and Weyl's lemma implies that it is a smooth harmonic function. Moreover, $v(x)=O(|x|^{-p})$ and $w(x)\to0$, so $w-v\to0$. Liouville's theorem therefore yields
\begin{equation}
w(x)=\int_{\partial F}|x-y|^{-p}\,d\eta(y),\qquad x\in\mathbb R^n.
\label{eq:2.3}
\end{equation}

It remains to determine the total mass of $\eta$. The kernel representation with compactly supported measure, together with its first derivatives, gives
$w(x)=O(|x|^{-p})$ and $\nabla_xw(x)=O(|x|^{-p-1})$. Green's formula on $B_R\setminus F$ gives
\[
\int_{B_R\setminus F}|\nabla_xw|^2\,dx
=-\int_{\partial F}\partial_\nu w\,dS
+\int_{\partial B_R}w\partial_rw\,dS.
\]
The absolute value of the last term is $O(R^{n-1-p-p-1})=O(R^{2-n})$, which tends to zero. Hence
\begin{equation}
\int_{\mathbb R^n}|\nabla_xw|^2\,dx
=-\int_{\partial F}\partial_\nu w\,dS
=\gamma_n\eta(\partial F).
\label{eq:2.4}
\end{equation}

To make explicit the relation between this Dirichlet minimum and the test class in Theorem~1.4, denote the energy in the preceding identity by $H$. For every admissible function $v$ in that theorem,
\[
\int_{\mathbb R^n}\nabla_xw\cdot\nabla_xv\,dx
=-\int_{\partial F}v\partial_\nu w\,dS\geq H.
\]
The Cauchy--Schwarz inequality gives
$H^2\leq H\int|\nabla_xv|^2$, and therefore
$\operatorname{cap}_{D,n}(F)\geq H$.
For the reverse inequality, let
$w_\delta=\min\{1,w/(1-\delta)\}$, where $0<\delta<1$.
This function is identically $1$ in a neighborhood of $F$, and
$\int|\nabla_xw_\delta|^2\leq H/(1-\delta)^2$.
Choose a smooth cutoff $\chi_R$ that equals $1$ on $B_R$, vanishes outside $B_{2R}$, and satisfies $|\nabla\chi_R|\leq C/R$. Since $w_\delta=O(|x|^{2-n})$,
\[
\int_{B_{2R}\setminus B_R}w_\delta^2|\nabla\chi_R|^2\,dx
=O(R^{2-n})\longrightarrow0.
\]
Together with square integrability of the gradient of $w_\delta$, this gives
$\int|\nabla(\chi_Rw_\delta)|^2\to\int|\nabla w_\delta|^2$.
For fixed $R$, convolve with a nonnegative smooth mollifier of sufficiently small support. The resulting functions remain identically $1$ in a neighborhood of $F$, belong to $C_c^\infty$, and their gradients converge in $L^2$. Letting first the mollification scale tend to zero, then $R\to\infty$, and finally $\delta\downarrow0$, we obtain
$\operatorname{cap}_{D,n}(F)\leq H$.
Consequently, Theorem~1.4 and \eqref{eq:2.4} give
\begin{equation}
\eta(\partial F)=\frac{\operatorname{cap}_{D,n}(F)}{\gamma_n}
=C_{n-2}(F)^{n-2}=c^p.
\label{eq:2.5}
\end{equation}

Now define the functions in one higher dimension by
\begin{equation}
u(X)=\int_{\partial F}R_y^{-p}\,d\eta(y),\qquad
\rho(X)=u(X)^{-1/p}.
\label{eq:2.6}
\end{equation}
By \eqref{eq:2.3}, $u(x,0)=w(x)$. Since $R_y\geq|x-y|$ and $\eta\geq0$,
$0<u(x,z)\leq w(x)\leq1$.
When $z\ne0$, every finite kernel value strictly decreases, and $\eta$ has positive total mass; hence $u(x,z)<w(x)$. When $z=0$ and $x\notin F$, we have $w(x)<1$. Thus $u=1$ on $F$, whereas $0<u<1$ on $\mathbb R^{n+1}\setminus F$. For $X_j\to X_0\in F$, Fatou's inequality for the nonnegative kernel and the global upper bound give
\[
1=u(X_0)\leq\liminf_j u(X_j)\leq\limsup_j u(X_j)\leq1.
\]
This proves continuity of $u$ at every point of $F$, and hence continuity of $\rho$ there.

\subsection{The exterior Green function in two dimensions}

When $n=2$, the inequality $c=C_0(F)>0$ ensures that $F$ has a logarithmic equilibrium measure $\sigma_0$ of finite energy. The Green function with pole at infinity can be defined in terms of the equilibrium potential by
\[
g(x)=\int_F\log|x-y|\,d\sigma_0(y)-\log c,\qquad x\in\Omega;
\]
see \cite[Definition~3.4 and equations~(3.1)--(3.2)]{ref7}.
This definition applies to any unbounded component of the complement and does not require $\Omega$ to be simply connected. Every point of a smooth boundary is regular for the Dirichlet problem, so $g$ extends continuously to the boundary with value $0$. Boundary regularity further implies that $g$ is smooth up to the boundary from the exterior. Since $\sigma_0$ is a compactly supported probability measure, expansion of the kernel at infinity gives
\begin{equation}
\Delta_xg=0\ \text{in }\Omega,\qquad
 g=0\ \text{on }\partial F,\qquad
 g(x)=\log|x|-\log c+O(|x|^{-1}).
\label{eq:2.7}
\end{equation}
Termwise differentiation of the same kernel also gives
\begin{equation}
\partial_rg=\frac1r+O(r^{-2}),\qquad r=|x|\longrightarrow\infty,
\label{eq:2.8}
\end{equation}
with an error uniform in the angular variable. On sufficiently large circles $g>0$. Applying the maximum principle in truncated exterior domains and then the strong maximum principle gives $g>0$ in $\Omega$. Set $g=0$ in $F$; then $g$ is continuous and nonnegative on the whole plane. The Hopf lemma gives $\partial_\nu g>0$ on every boundary component.

Define the boundary flux measure
\begin{equation}
d\mu=\frac{\partial_\nu g}{2\pi}\,dS.
\label{eq:2.9}
\end{equation}
This is a finite nonnegative measure with smooth density. For every $\varphi\in C_c^\infty(\mathbb R^2)$, integration by parts gives
\[
\langle\Delta_xg,\varphi\rangle
=-\int_\Omega\nabla_xg\cdot\nabla_x\varphi\,dx
=\int_{\partial F}\varphi\partial_\nu g\,dS
=2\pi\int\varphi\,d\mu.
\]
On the other hand, the divergence theorem on $B_R\setminus F$ gives
\[
0=\int_{B_R\setminus F}\Delta_xg\,dx
=\int_{\partial B_R}\partial_rg\,dS
-\int_{\partial F}\partial_\nu g\,dS.
\]
By \eqref{eq:2.8}, the integral on the outer circle is $2\pi+O(R^{-1})$. Letting $R\to\infty$ yields
\begin{equation}
\mu(\partial F)=\frac1{2\pi}\int_{\partial F}\partial_\nu g\,dS=1.
\label{eq:2.10}
\end{equation}

Set $v(x)=\int_{\partial F}\log|x-y|\,d\mu(y)$.
The density on the smooth curve is bounded, and the length of the curve in a ball of radius $r$ is at most $Cr$. For $0<t<1/2$, the same annular decomposition gives the uniform estimate
\[
\begin{aligned}
\int_{\partial F\cap B(x,t)}\bigl|\log|x-y|\bigr|\,d\mu(y)
&\leq Ct\sum_{k=0}^\infty2^{-k}\bigl(|\log t|+(k+1)\log2\bigr)\\
&\leq C't(1+|\log t|).
\end{aligned}
\]
Fix $x_0$, and take $0<t<1/6$ and $|x-x_0|<t$. Then the fixed near region satisfies
$B(x_0,2t)\subset B(x,3t)$, and both near-region integrals are
$O(t(1+|\log t|))$. On its complement the kernel stays a positive distance from its singularity, so the dominated convergence theorem applies. Letting first $x\to x_0$ and then $t\downarrow0$ shows that $v$ is finite everywhere and continuous. The identity
$\Delta_x\log|x-y|=2\pi\delta_y$ gives
$\Delta_xv=2\pi\mu=\Delta_xg$; thus $v-g$ is distributionally harmonic on the whole plane. Furthermore, \eqref{eq:2.10} and compact support give
$v(x)=\log|x|+O(|x|^{-1})$, so $v(x)-g(x)\to\log c$.
Weyl's lemma and Liouville's theorem imply $v-g\equiv\log c$. In other words,
\begin{equation}
g(x)=\int_{\partial F}\log|x-y|\,d\mu(y)-\log c,
\qquad x\in\mathbb R^2.
\label{eq:2.11}
\end{equation}
All the integrals, distributional identities, and divergence formulas above include all boundary components, so they apply to multiply connected exterior domains.

Define
\begin{equation}
G(X)=\int_{\partial F}\log R_y\,d\mu(y)-\log c,\qquad
\rho(X)=e^{G(X)}.
\label{eq:2.12}
\end{equation}
Since $\log R_y\geq\log|x-y|$, \eqref{eq:2.11} implies
$G(x,z)\geq g(x)\geq0$. The integral strictly increases when $z\ne0$; when $z=0$ and $x\notin F$, we have $g(x)>0$.
Hence $G=0$ on $F$, whereas $G>0$ on $\mathbb R^3\setminus F$.
For $X_j=(x_j,z_j)\to(x_0,0)\in F$, the inequality
$G(X_j)\geq g(x_j)\to0$ gives $\liminf_jG(X_j)\geq0$.
For these $X_j$ and $y$ in the compact set $\partial F$, the functions $\log R_y$ have a common finite upper bound $M$. Applying Fatou's inequality to the nonnegative functions $M-\log R_y$ gives
\[
\limsup_j\int\log|X_j-(y,0)|\,d\mu(y)
\leq\int\log|x_0-y|\,d\mu(y)=\log c.
\]
It follows that $\limsup_jG(X_j)\leq0$, which proves continuity of $G$ and $\rho=e^G$ at every point of $F$.

\subsection{A common differential equation}

For any $X_0\notin F$, some neighborhood of $X_0$ has positive distance from the compact set $\partial F\times\{0\}$. All derivatives of the kernels are bounded on the product of this neighborhood and the support of integration, so differentiation under the integral sign is valid to every order. Thus $u$ or $G$, and also $\rho$, are smooth on $\mathbb R^{n+1}\setminus F$. Since $R_y^2=|x-y|^2+z^2$, all these functions are even in $z$.

For $z>0$, introduce the elliptic operator whose drift coefficient is singular at $z=0$:
\begin{equation}
\mathcal E=\Delta_x+\partial_z^2-\frac1z\partial_z
=\Delta-\frac1z\partial_z.
\label{eq:2.13}
\end{equation}
When $n\geq3$ and $p=n-2$, direct calculation gives
\[
\Delta R_y^{-p}
=p\bigl(p+2-(n+1)\bigr)R_y^{-p-2}=-pR_y^{-p-2},
\qquad
\frac1z\partial_zR_y^{-p}=-pR_y^{-p-2}.
\]
Hence $\mathcal E u=0$. When $n=2$, the extended ambient space is $\mathbb R^3$, and
\[
\Delta\log R_y=R_y^{-2},\qquad
\frac1z\partial_z\log R_y=R_y^{-2}.
\]
Hence $\mathcal E G=0$.

\begin{proposition}[Basic properties of the extension]
\label{prop:extension}
The function $\rho$ defined by \eqref{eq:2.6} or \eqref{eq:2.12} is continuous on the whole space, even in $z$, and smooth outside $F$. It satisfies
\begin{equation}
\rho=1\ \text{on }F,\qquad
\rho>1\ \text{outside }F,\qquad
\rho(X)=\frac{|X|}{c}+O(1),\quad |X|\longrightarrow\infty.
\label{eq:2.14}
\end{equation}
For $z\ne0$,
\begin{equation}
\Delta\rho=\frac{n-1}{\rho}|\nabla\rho|^2+\frac{\rho_z}{z}.
\label{eq:2.15}
\end{equation}
At $z=0$ and $x\notin F$, this identity extends continuously, with $\rho_z/z$ assigned the value $\rho_{zz}(x,0)$.
\end{proposition}

\begin{proof}
The range, continuity, smoothness, and evenness have already been proved. We verify the asymptotic formula at infinity and the differential equation. Suppose that the support of integration is contained in $\{|y|\leq M\}$, and let $r=|X|>2M$. Since
$\bigl||X-(y,0)|-r\bigr|\leq M$, uniformly on the support of integration,
\[
R_y^{-p}=r^{-p}\bigl(1+O(r^{-1})\bigr),\qquad
\log R_y=\log r+O(r^{-1}).
\]
When $n\geq3$, \eqref{eq:2.5} gives
$u(X)=c^pr^{-p}(1+O(r^{-1}))$. Taking the power $-1/p$ yields
$\rho(X)=(r/c)(1+O(r^{-1}))=r/c+O(1)$.
When $n=2$, \eqref{eq:2.10} gives $G(X)=\log(r/c)+O(r^{-1})$.
Exponentiation again yields $\rho(X)=r/c+O(1)$.

When $n\geq3$, substitute $u=\rho^{-p}$ for $z>0$ and apply the chain rule term by term to obtain
\[
0=\mathcal E(\rho^{-p})
=-p\rho^{-p-1}\mathcal E\rho
+p(p+1)\rho^{-p-2}|\nabla\rho|^2.
\]
Since $p>0$, $\rho>0$, and $p+1=n-1$, division by $p\rho^{-p-2}$ gives
$\rho\mathcal E\rho=(n-1)|\nabla\rho|^2$.
When $n=2$, we have $G=\log\rho$, and hence
\[
0=\mathcal E(\log\rho)
=\rho^{-1}\mathcal E\rho-\rho^{-2}|\nabla\rho|^2.
\]
This gives the same identity, since now $n-1=1$.
Expanding $\mathcal E\rho$ yields \eqref{eq:2.15}.
Evenness in $z$ shows that the same formula holds for $z<0$.

Finally, fix $x\notin F$. Smoothness and evenness give $\rho_z(x,0)=0$, and
\[
\frac{\rho_z(x,z)}{z}
=\int_0^1\rho_{zz}(x,tz)\,dt
\longrightarrow\rho_{zz}(x,0).
\]
All other terms are also continuous, so \eqref{eq:2.15} extends to this point in the stated manner.
\end{proof}

\section{A pointwise gradient estimate}

In this section we use bounds for the potentials in the whole space to prove the gradient inequality needed for the construction of a superharmonic function. Here the gradient is taken with respect to all variables $X=(x,z)\in\mathbb R^{n+1}$.

\begin{proposition}[Pointwise gradient estimate]
Let $\rho$ be the function constructed in Section~2. For every $X=(x,z)$ with $z\ne0$, we have
\begin{equation}
\rho|\nabla\rho|^2\le(\rho^2-1)\frac{\rho_z}{z}.
\label{eq:3.1}
\end{equation}
The inequality also holds at points with $z=0$ and $x\notin F$, where $\rho_z/z$ is assigned its continuous extension $\rho_{zz}(x,0)$.
\end{proposition}

\begin{proof}
First suppose that $n\ge3$, and fix $z\ne0$. Put $p=n-2$, retain the measure $\eta$ from Section~2, and set
\begin{equation}
\begin{aligned}
S&=\int_{\partial F}R_y^{-p-2}\,d\eta(y),\\
M&=\int_{\partial F}(X-(y,0))R_y^{-p-2}\,d\eta(y)=(M_h,zS),
\end{aligned}
\label{eq:3.2}
\end{equation}
where $M_h\in\mathbb R^n$ is the horizontal component of $M$. Since $R_y\ge|z|>0$ and $\eta$ is a compactly supported, finite, nonzero positive measure, these integrals are finite and $S>0$.

Differentiating the integral representation shows that $M$ determines the full gradient of $u$, whereas $S$ determines its vertical derivative. We therefore seek to bound $|M|^2$ in terms of $u$ and $S$. To this end, we introduce a quadratic integral depending on a horizontal point $b$. Completing the square will show that its minimum is $u-|M|^2/S$, while the bound for the capacitary potential will give a lower bound valid for every $b$.

More precisely, for any $b\in\mathbb R^n$, define
\[
J_b=\int_{\partial F}|b-y|^2R_y^{-p-2}\,d\eta(y).
\]
The flux representation and the maximum principle for the exterior capacitary potential give
\[
\int_{\partial F}|b-y|^{-p}\,d\eta(y)=w(b)\le1.
\]
This identity is also valid when $b\in\partial F$. Indeed, $\eta$ has a bounded density with respect to the surface measure on the smooth boundary, and $p=n-2$. Thus, in local coordinates on the boundary, the integral near the singularity is bounded by
\[
C\int_0^\varepsilon r^{n-2}r^{-p}\,dr=C\varepsilon.
\]
In particular, the kernel is locally integrable, and each point of the boundary has $\eta$-measure zero.

For $y\ne b$, we have
\[
R_y^{-p}
=\bigl(|b-y|^{-p}\bigr)^{2/(p+2)}
\bigl(|b-y|^2R_y^{-p-2}\bigr)^{p/(p+2)}.
\]
Applying H\"older's inequality with conjugate exponents $(p+2)/2$ and $(p+2)/p$, we obtain
\begin{equation}
\begin{aligned}
u(X)
&\le\left(\int_{\partial F}|b-y|^{-p}\,d\eta(y)\right)^{2/(p+2)}
\left(\int_{\partial F}|b-y|^2R_y^{-p-2}\,d\eta(y)\right)^{p/(p+2)}\\
&=w(b)^{2/(p+2)}J_b^{p/(p+2)}\le J_b^{p/(p+2)}.
\end{aligned}
\label{eq:3.3}
\end{equation}
Consequently, for every $b\in\mathbb R^n$,
\[
J_b\ge u^{(p+2)/p}=u^{1+2/p}.
\]

Let $a=x-b$. Since $b-y=(x-y)-a$ and
\[
\int_{\partial F}|x-y|^2R_y^{-p-2}\,d\eta(y)
=\int_{\partial F}(R_y^2-z^2)R_y^{-p-2}\,d\eta(y)
=u-z^2S,
\]
expanding and then completing the square yields
\begin{equation}
\begin{aligned}
J_b&=u-z^2S-2a\cdot M_h+S|a|^2\\
&=u-z^2S-\frac{|M_h|^2}{S}
+S\left|a-\frac{M_h}{S}\right|^2\\
&=u-\frac{|M|^2}{S}
+S\left|a-\frac{M_h}{S}\right|^2.
\end{aligned}
\label{eq:3.4}
\end{equation}
In the last step we used $|M|^2=|M_h|^2+z^2S^2$. Since $b$ ranges over all of $\mathbb R^n$, we may choose $a=M_h/S$. Combining this choice with $J_b\ge u^{1+2/p}$ gives
\[
u-\frac{|M|^2}{S}\ge u^{1+2/p},
\]
or equivalently,
\begin{equation}
|M|^2\le uS(1-u^{2/p}).
\label{eq:3.5}
\end{equation}

Our fixed point $X$ satisfies $z\ne0$, so the integral kernels and their derivatives have a common integrable majorant in a neighborhood of $X$. Differentiation under the integral sign is therefore justified and gives
\[
\nabla u=-pM,\qquad u_z=-pzS.
\]
Since $\rho=u^{-1/p}$, the chain rule gives
\[
\nabla\rho=-\frac1p u^{-1/p-1}\nabla u=\frac\rho u M,
\qquad
\frac{\rho_z}{z}=\frac\rho u S.
\]
Substituting \eqref{eq:3.5} and using $u^{2/p}=\rho^{-2}$, we obtain
\[
\begin{aligned}
\rho|\nabla\rho|^2
&=\frac{\rho^3}{u^2}|M|^2\\
&\le\frac{\rho^3S}{u}(1-u^{2/p})
=\frac{\rho^3S}{u}(1-\rho^{-2})\\
&=(\rho^2-1)\frac{\rho S}{u}
=(\rho^2-1)\frac{\rho_z}{z}.
\end{aligned}
\]

Now suppose that $n=2$, again with $z\ne0$ fixed. Retain the probability measure $\mu$ in \eqref{eq:2.9} and \eqref{eq:2.10}, and set
\[
S=\int_{\partial F}R_y^{-2}\,d\mu(y),
\qquad
M=\int_{\partial F}(X-(y,0))R_y^{-2}\,d\mu(y)=(M_h,zS).
\]
Again $S>0$, and these integrals are finite. We still need to control $|M|^2$ by $S$. The corresponding quadratic integral has minimum $1-|M|^2/S$, and a lower bound follows from the logarithmic potential and Jensen's inequality. For any $b\in\mathbb R^2$, put
\[
J_b=\int_{\partial F}\frac{|b-y|^2}{R_y^2}\,d\mu(y).
\]
The measure $\mu$ has a bounded density with respect to arclength on the smooth curve, and
\[
\int_0^\varepsilon|\log r|\,dr<\infty.
\]
Hence $\log|b-y|$ is $\mu$-integrable even when $b\in\partial F$. Also, $\log R_y$ is bounded because $R_y\ge|z|$ and $\partial F$ is compact. We may therefore apply Jensen's inequality for the logarithm with respect to the probability measure $\mu$, obtaining
\begin{equation}
\begin{aligned}
J_b
&\ge\exp\left(\int_{\partial F}\log\frac{|b-y|^2}{R_y^2}\,d\mu(y)\right)\\
&=\exp\left(2\int_{\partial F}\log|b-y|\,d\mu(y)
-2\int_{\partial F}\log R_y\,d\mu(y)\right)\\
&=\exp\bigl(2(g(b)+\log c)-2(G(X)+\log c)\bigr)\\
&=e^{2g(b)-2G(X)}\ge e^{-2G(X)}=\rho^{-2}.
\end{aligned}
\label{eq:3.6}
\end{equation}
Here we used the inequality $g(b)\ge0$ on the whole plane.

Again let $a=x-b$. Using $\mu(\partial F)=1$, we first compute
\[
\int_{\partial F}|x-y|^2R_y^{-2}\,d\mu(y)
=\int_{\partial F}(1-z^2R_y^{-2})\,d\mu(y)
=1-z^2S.
\]
It follows that
\[
\begin{aligned}
J_b&=1-z^2S-2a\cdot M_h+S|a|^2\\
&=1-z^2S-\frac{|M_h|^2}{S}
+S\left|a-\frac{M_h}{S}\right|^2\\
&=1-\frac{|M|^2}{S}
+S\left|a-\frac{M_h}{S}\right|^2.
\end{aligned}
\]
Taking $a=M_h/S$ in \eqref{eq:3.6}, we obtain
\[
1-\frac{|M|^2}{S}\ge\rho^{-2},
\qquad
|M|^2\le S(1-\rho^{-2}).
\]
Differentiating the integral representation of $G$ gives $\nabla G=M$ and $G_z=zS$. Since $\rho=e^G$,
\[
\nabla\rho=\rho M,\qquad\frac{\rho_z}{z}=\rho S.
\]
Therefore,
\[
\begin{aligned}
\rho|\nabla\rho|^2
&=\rho^3|M|^2\le\rho^3S(1-\rho^{-2})\\
&=(\rho^2-1)\rho S
=(\rho^2-1)\frac{\rho_z}{z}.
\end{aligned}
\]

Finally, if $z=0$ and $x\notin F$, then $\operatorname{dist}(x,F)>0$. Thus the kernels above are smooth in a neighborhood of $(x,0)$ and may be differentiated under the integral sign. Since $\rho$ is even in $z$, we have $\rho_z(x,0)=0$ and
\[
\frac{\rho_z(x,z)}{z}
=\int_0^1\rho_{zz}(x,tz)\,dt
\longrightarrow\rho_{zz}(x,0).
\]
Letting $z\to0$ in the inequality already proved gives the asserted inequality on this plane.
\end{proof}

\section{An ordinary differential equation and a supersolution}

For $\Phi\in C^2((1,\infty))$, the chain rule and \eqref{eq:2.15} give, outside $F$,
\begin{equation}
\begin{aligned}
\Delta\Phi(\rho)
&=\Phi''(\rho)|\nabla\rho|^2+\Phi'(\rho)\Delta\rho\\
&=\left(\Phi''(\rho)+\frac{n-1}{\rho}\Phi'(\rho)\right)|\nabla\rho|^2
+\Phi'(\rho)\frac{\rho_z}{z}.
\end{aligned}
\label{eq:4.1}
\end{equation}
At points with $z=0$ and $x\notin F$, the last term is again interpreted by continuous extension.

By the gradient estimate, outside $F$ we have
\[
\frac{\rho_z}{z}-\frac{\rho|\nabla\rho|^2}{\rho^2-1}\geq0.
\]
To construct a decreasing superharmonic function $U=\Phi_n(\rho)$, we choose $\Phi_n$ so that the right-hand side of \eqref{eq:4.1} is the product of this nonnegative quantity and $\Phi_n'(\rho)<0$.
Comparing the coefficients of $|\nabla\rho|^2$ and imposing the boundary values on the conductor and at infinity, we obtain
\begin{equation}
\begin{gathered}
\Phi_n''(t)+\left(\frac{n-1}{t}+\frac{t}{t^2-1}\right)\Phi_n'(t)=0
\quad(t>1),\\
\Phi_n(1)=1,\qquad \lim_{t\to\infty}\Phi_n(t)=0.
\end{gathered}
\label{eq:4.2}
\end{equation}
Here the endpoint value $\Phi_n(1)$ is defined by continuous extension.

\begin{lemma}[Solution of the model equation]\label{lem:4.1}
Let
\begin{equation}
Z_n=\int_1^\infty\frac{t^{1-n}}{\sqrt{t^2-1}}\,dt
=\frac{\sqrt{\pi}\,\Gamma((n-1)/2)}{2\Gamma(n/2)}.
\label{eq:4.3}
\end{equation}
The unique solution of the boundary-value problem \eqref{eq:4.2} in
$C^2((1,\infty))\cap C([1,\infty))$ is
\begin{equation}
\Phi_n(r)=\frac{1}{Z_n}\int_r^\infty\frac{t^{1-n}}{\sqrt{t^2-1}}\,dt,
\qquad r\geq1.
\label{eq:4.4}
\end{equation}
It is strictly decreasing on $(1,\infty)$, and
\begin{equation}
\Phi_n(r)=a_n r^{-(n-1)}\bigl(1+O(r^{-2})\bigr),\qquad
a_n=\frac{1}{(n-1)Z_n}
=\frac{\Gamma(n/2)}{\sqrt{\pi}\,\Gamma((n+1)/2)}.
\label{eq:4.5}
\end{equation}
Moreover, direct differentiation with respect to the lower limit of integration yields
\begin{equation}
\Phi_n'(r)=-\frac{1}{Z_n}\frac{r^{1-n}}{\sqrt{r^2-1}}
=-(n-1)a_n r^{-n}\bigl(1+O(r^{-2})\bigr).
\label{eq:4.6}
\end{equation}
\end{lemma}

\begin{proof}
Equation \eqref{eq:4.2} is equivalent to
\[
\frac{d}{dt}\left(t^{n-1}\sqrt{t^2-1}\,\Phi_n'(t)\right)=0.
\]
Hence there is a constant $C$ such that
\[
\Phi_n'(t)=C\frac{t^{1-n}}{\sqrt{t^2-1}}.
\]
The integrand has order $(t-1)^{-1/2}$ near $t=1$ and order $t^{-n}$ at infinity.
Since $n\geq2$, it is integrable at both endpoints.
Integrating and using the boundary value at infinity gives
\[
\Phi_n(r)=-C\int_r^\infty\frac{t^{1-n}}{\sqrt{t^2-1}}\,dt.
\]
Using $\Phi_n(1)=1$ then gives $C=-1/Z_n$, proving \eqref{eq:4.4} and uniqueness.
Since $Z_n>0$, the derivative is strictly negative.

To evaluate $Z_n$, set $s=t^{-2}$. Then
\[
t=s^{-1/2},\qquad dt=-\tfrac12s^{-3/2}\,ds,\qquad
\sqrt{t^2-1}=s^{-1/2}\sqrt{1-s}.
\]
Thus
\[
Z_n=\frac12\int_0^1s^{(n-3)/2}(1-s)^{-1/2}\,ds.
\]
For $\alpha,\beta>0$, write
\[
B(\alpha,\beta)=\int_0^1s^{\alpha-1}(1-s)^{\beta-1}\,ds.
\]
In the double integral for $\Gamma(\alpha)\Gamma(\beta)$, make the change of variables
$v=\xi+\zeta$, $s=\xi/(\xi+\zeta)$, whose Jacobian has absolute value $v$.
This gives
\[
\begin{aligned}
\Gamma(\alpha)\Gamma(\beta)
&=\int_0^\infty\int_0^\infty
e^{-(\xi+\zeta)}\xi^{\alpha-1}\zeta^{\beta-1}\,d\xi\,d\zeta\\
&=\left(\int_0^\infty e^{-v}v^{\alpha+\beta-1}\,dv\right)
\left(\int_0^1s^{\alpha-1}(1-s)^{\beta-1}\,ds\right)\\
&=\Gamma(\alpha+\beta)B(\alpha,\beta).
\end{aligned}
\]
Taking $\alpha=(n-1)/2$ and $\beta=1/2$, and using $\Gamma(1/2)=\sqrt{\pi}$, we obtain
\[
Z_n=\frac12B\left(\frac{n-1}{2},\frac12\right)
=\frac{\sqrt{\pi}\,\Gamma((n-1)/2)}{2\Gamma(n/2)}.
\]

As $t\to\infty$,
\[
\frac{t^{1-n}}{\sqrt{t^2-1}}
=t^{-n}(1-t^{-2})^{-1/2}
=t^{-n}+\frac12t^{-n-2}+O(t^{-n-4}).
\]
Integrating term by term and estimating the absolute value of the remainder gives
\[
\Phi_n(r)=\frac{r^{1-n}}{(n-1)Z_n}
+\frac{r^{-n-1}}{2(n+1)Z_n}+O(r^{-n-3}).
\]
This proves the asymptotic expansion in \eqref{eq:4.5}.
The expression for $a_n$ follows from
$\Gamma((n+1)/2)=((n-1)/2)\Gamma((n-1)/2)$.
Finally, the fundamental theorem of calculus gives the exact derivative formula in \eqref{eq:4.6}; expanding $(1-r^{-2})^{-1/2}$ then gives the stated asymptotic expansion of the derivative.
\end{proof}

Set $U=\Phi_n(\rho)$. By \eqref{eq:4.2},
\[
\Phi_n''(\rho)+\frac{n-1}{\rho}\Phi_n'(\rho)
=-\frac{\rho}{\rho^2-1}\Phi_n'(\rho).
\]
Substituting this into \eqref{eq:4.1} yields
\begin{equation}
\Delta U=\Phi_n'(\rho)
\left(\frac{\rho_z}{z}-\frac{\rho|\nabla\rho|^2}{\rho^2-1}\right)\leq0.
\label{eq:4.7}
\end{equation}
Indeed, $\rho>1$ outside $F$, so \eqref{eq:3.1} shows that the quantity in parentheses is nonnegative, while Lemma~\ref{lem:4.1} gives $\Phi_n'(\rho)<0$.
Thus $U$ is superharmonic outside $F$, that is, $-\Delta U\geq0$.
Since $\rho$ is continuous, $\rho|_F=1$, and $\Phi_n(1)=1$, the function $U$ is continuous throughout the whole space and satisfies
\[
U=1\quad\text{on }F,\qquad
0<U<1\quad\text{on }\mathbb{R}^{n+1}\setminus F.
\]

We next derive the asymptotic expansions of $U$ and its radial derivative at infinity separately from the flux integral representation.
Write $X=r\theta$, where $r=|X|$ and $\theta\in S^n$, and let $\theta_h\in\mathbb{R}^n$ denote the horizontal component of $\theta$.
Below, $\partial_r$ denotes differentiation with respect to $r$ while keeping $\theta$ fixed.
Since the boundary is compact, we may choose $L>0$ such that $|y|\leq L$ for every $y\in\partial F$.
All expansions are taken for $r>2L$, and the constants in the remainder terms are independent of $\theta$.

Suppose $n\geq3$, set $p=n-2$, and define the barycenter of the flux measure by
\[
m=\frac{1}{c^p}\int_{\partial F}y\,d\eta(y)\in\mathbb{R}^n.
\]
Writing $k=\theta_h\cdot y$, we have
\[
R_y^2=r^2-2rk+|y|^2.
\]
Taylor expansion of $(1+s)^{-p/2}$ gives
\[
R_y^{-p}=r^{-p}\left(1+\frac{pk}{r}+O(r^{-2})\right).
\]
On the other hand, direct differentiation of the kernel gives
\[
\begin{aligned}
\partial_r(R_y^{-p})
&=-p(r-k)(r^2-2rk+|y|^2)^{-(p+2)/2}\\
&=-pr^{-p-1}\left(1-\frac{k}{r}\right)
\left(1+\frac{(p+2)k}{r}+O(r^{-2})\right)\\
&=-pr^{-p-1}\left(1+\frac{(p+1)k}{r}+O(r^{-2})\right).
\end{aligned}
\]
Integrating these expressions separately and using $\eta(\partial F)=c^p$, we obtain
\[
\begin{aligned}
u(r\theta)&=c^pr^{-p}\left(1+\frac{p\theta_h\cdot m}{r}+O(r^{-2})\right),\\
\partial_r u(r\theta)&=-pc^pr^{-p-1}
\left(1+\frac{(p+1)\theta_h\cdot m}{r}+O(r^{-2})\right).
\end{aligned}
\]
Therefore
\[
\begin{aligned}
\rho(r\theta)
&=\frac{r}{c}\left(1+\frac{p\theta_h\cdot m}{r}+O(r^{-2})\right)^{-1/p}\\
&=\frac{r}{c}-\frac{\theta_h\cdot m}{c}+O(r^{-1}),\\
\partial_r\rho(r\theta)
&=-\frac1p u^{-1/p-1}\partial_r u\\
&=\frac1c\left(1-\frac{(p+1)\theta_h\cdot m}{r}+O(r^{-2})\right)
\left(1+\frac{(p+1)\theta_h\cdot m}{r}+O(r^{-2})\right)\\
&=\frac1c+O(r^{-2}).
\end{aligned}
\]

When $n=2$, set $m=\int_{\partial F}y\,d\mu(y)$.
Starting from the same identity $R_y^2=r^2-2rk+|y|^2$, we obtain, respectively,
\[
\begin{aligned}
\log R_y
&=\log r+\frac12\log\left(1-\frac{2k}{r}+\frac{|y|^2}{r^2}\right)
=\log r-\frac{k}{r}+O(r^{-2}),\\
\partial_r\log R_y
&=\frac{r-k}{r^2-2rk+|y|^2}
=\frac1r+\frac{k}{r^2}+O(r^{-3}).
\end{aligned}
\]
Since $\mu(\partial F)=1$, integration gives
\[
G(r\theta)=\log(r/c)-\frac{\theta_h\cdot m}{r}+O(r^{-2}),\qquad
\partial_r G(r\theta)=\frac1r+\frac{\theta_h\cdot m}{r^2}+O(r^{-3}).
\]
As $\rho=e^G$,
\[
\begin{aligned}
\rho(r\theta)
&=\frac{r}{c}\left(1-\frac{\theta_h\cdot m}{r}+O(r^{-2})\right)
=\frac{r}{c}-\frac{\theta_h\cdot m}{c}+O(r^{-1}),\\
\partial_r\rho(r\theta)
&=\rho\,\partial_r G\\
&=\frac1c\left(1-\frac{\theta_h\cdot m}{r}+O(r^{-2})\right)
\left(1+\frac{\theta_h\cdot m}{r}+O(r^{-2})\right)
=\frac1c+O(r^{-2}).
\end{aligned}
\]

In all dimensions $n\geq2$, we have therefore shown that
\[
\rho(r\theta)=\frac{r}{c}\bigl(1+O(r^{-1})\bigr),\qquad
\partial_r\rho(r\theta)=\frac1c+O(r^{-2}),
\]
with remainder terms uniform in the direction.
Equations \eqref{eq:4.5} and \eqref{eq:4.6}, together with the chain rule, now give
\[
\begin{aligned}
U(r\theta)
&=a_n\rho^{-(n-1)}\bigl(1+O(\rho^{-2})\bigr)
=a_nc^{n-1}r^{-(n-1)}\bigl(1+O(r^{-1})\bigr),\\
\partial_r U(r\theta)
&=\Phi_n'(\rho)\partial_r\rho\\
&=-(n-1)a_n\rho^{-n}\bigl(1+O(\rho^{-2})\bigr)
\left(\frac1c+O(r^{-2})\right)\\
&=-(n-1)a_nc^{n-1}r^{-n}\bigl(1+O(r^{-1})\bigr).
\end{aligned}
\]
That is,
\begin{equation}
\begin{aligned}
U(X)&=a_nc^{n-1}|X|^{-(n-1)}+O(|X|^{-n}),\\
\partial_r U(X)&=-(n-1)a_nc^{n-1}|X|^{-n}+O(|X|^{-n-1}).
\end{aligned}
\label{eq:4.8}
\end{equation}

\subsection{Explicit formulas in dimensions two and three}

When $n=2$, we have $Z_2=\pi/2$.
Making the substitution $s=1/t$ in \eqref{eq:4.4}, we obtain
\[
\Phi_2(\rho)=\frac2\pi\int_0^{1/\rho}\frac{ds}{\sqrt{1-s^2}}
=\frac2\pi\arcsin(\rho^{-1}),\qquad a_2=\frac2\pi.
\]
When $n=3$, we have $Z_3=1$.
The same substitution gives
\[
\Phi_3(\rho)=\int_0^{1/\rho}\frac{s\,ds}{\sqrt{1-s^2}}
=1-\sqrt{1-\rho^{-2}},\qquad a_3=\frac12.
\]
In this case $p=1$ and $\rho^{-1}=u$, so $U=1-\sqrt{1-u^2}$.
In Section~5, the coefficient in \eqref{eq:4.8} will yield
$C_2(K)^2\leq C_1(K)^2/2$.

\section{The Dirichlet energy estimate and proof of the main theorem}

We first establish flux and energy identities on bounded truncations and then construct smooth test functions admissible in the definition of capacity. We next prove smooth outer approximation and continuity of capacity from above, and complete the proof of the main theorem.

\subsection{An energy bound for supersolutions}

\begin{lemma}[An energy bound for supersolutions]
\label{lem:5.1}
Let $E\subset\mathbb R^d$ be compact, where $d\geq3$. Suppose that $U$ is continuous on $\mathbb R^d$, smooth on $\mathbb R^d\setminus E$, and satisfies
\[
 U=1\quad\text{on }E,\qquad 0<U\leq1,\qquad
 -\Delta U\geq0\quad\text{on }\mathbb R^d\setminus E.
\]
Write $r=|X|$. If, for some constant $a>0$, the expansions
\begin{equation}
 U(X)=ar^{2-d}+O(r^{1-d}),\qquad
 \partial_rU(X)=-(d-2)ar^{1-d}+O(r^{-d})
 \label{eq:5.1}
\end{equation}
hold uniformly in direction, then
\begin{equation}
 \operatorname{cap}_{D,d}(E)\leq(d-2)|S^{d-1}|a.
 \label{eq:5.2}
\end{equation}
\end{lemma}

\begin{proof}
The conclusion is immediate when $E=\varnothing$, so assume that $E\neq\varnothing$. Set
\[
 M_R=-\int_{\partial B_R}\partial_rU\,dS.
\]
The second expansion in \eqref{eq:5.1}, with its remainder uniform in direction, gives
\begin{equation}
 M_R=(d-2)|S^{d-1}|a+O(R^{-1})
 \longrightarrow M:=(d-2)|S^{d-1}|a.
 \label{eq:5.3}
\end{equation}

Choose a regular value $t\in(0,1)$ of the smooth function $U|_{\mathbb R^d\setminus E}$; that is, $\nabla U(X)\neq0$ at every point where $U(X)=t$. Define
\[
 D_t=\{X:U(X)<t\},\qquad S_t=\{X:U(X)=t\},\qquad
 J_t=\int_{S_t}|\nabla U|\,dS.
\]
Continuity and the identity $U=1$ on $E$ ensure that $E$ has an open neighborhood whose closure is disjoint from $D_t$. Moreover, since $U\to0$ at infinity, the set $\{U\geq t\}$ is bounded. Thus $S_t$ is a compact smooth hypersurface contained in $\mathbb R^d\setminus E$, and $J_t<\infty$. A compact smooth manifold can be covered by finitely many connected coordinate neighborhoods, so $S_t$ has only finitely many connected components.

Choose $R_t>0$ such that $\{U\geq t\}\subset B_{R_t}$. For $R>R_t$, set $\Omega_{t,R}=D_t\cap B_R$. Its boundary is the disjoint union of $S_t$ and $\partial B_R$. On $S_t$, the outward normal relative to $\Omega_{t,R}$ is $\nabla U/|\nabla U|$, and hence
\[
 \partial_\nu U=|\nabla U|\quad\text{on }S_t.
\]
On $\partial B_R$, we have $\partial_\nu U=\partial_rU$. If the truncated domain has several connected components, the divergence formulas below are applied to each component and then summed, with all boundary components included.

Set $f=-\Delta U\geq0$. Since $\overline{\Omega_{t,R}}$ is disjoint from $E$, the function $U$ is smooth in a neighborhood of this closed set. The divergence theorem therefore gives the identity on the bounded domain
\[
 -\int_{\Omega_{t,R}}f\,dX
 =\int_{\Omega_{t,R}}\Delta U\,dX
 =J_t+\int_{\partial B_R}\partial_rU\,dS
 =J_t-M_R.
\]
Let $R\to\infty$. The corresponding integral of the nonnegative function $f$ increases monotonically with $R$, and the right-hand side has a finite limit. The monotone convergence theorem yields
\begin{equation}
 \int_{D_t}f\,dX=M-J_t<\infty.
 \label{eq:5.4}
\end{equation}
In particular, $0\leq J_t\leq M$. Since $0<U<t$ on $D_t$, we also have $Uf\in L^1(D_t)$. No global integrability of $\nabla U$ has been used at this point.

On the same bounded domain, the identity
\[
 \operatorname{div}(U\nabla U)=|\nabla U|^2+U\Delta U
 =|\nabla U|^2-Uf
\]
gives
\[
 \int_{\Omega_{t,R}}|\nabla U|^2\,dX
 =tJ_t+\int_{\Omega_{t,R}}Uf\,dX
   +\int_{\partial B_R}U\partial_rU\,dS.
\]
By the two expansions in \eqref{eq:5.1}, the absolute value of the last term is $O(R^{2-d})$, which tends to zero because $d\geq3$. The second term is already known to be integrable, and the first term is finite. Another application of the monotone convergence theorem establishes both $\nabla U\in L^2(D_t)$ and
\begin{equation}
 \int_{D_t}|\nabla U|^2\,dX
 =tJ_t+\int_{D_t}Uf\,dX.
 \label{eq:5.5}
\end{equation}
Substituting \eqref{eq:5.4}, we obtain
\begin{equation}
 \int_{D_t}|\nabla U|^2\,dX
 =tM-\int_{D_t}(t-U)f\,dX\leq tM.
 \label{eq:5.6}
\end{equation}

Define
\[
 v_t=\min\{U/t,1\}.
\]
The function $v_t$ is identically $1$ in a neighborhood of $E$; near every other point, it is a Lipschitz truncation of a smooth function. Thus $v_t$ is locally Lipschitz on $\mathbb R^d$, and in particular $v_t\in W^{1,2}_{\mathrm{loc}}(\mathbb R^d)$. Here $W^{1,2}$ denotes the Sobolev space of square-integrable functions with square-integrable first weak derivatives, and the subscript $\mathrm{loc}$ means that this condition holds locally. We write $\mathbf1_A$ for the indicator of a set $A$, equal to $1$ on $A$ and $0$ on its complement. Since the regular level set $S_t$ has zero $d$-dimensional Lebesgue measure, the Sobolev chain rule gives the weak gradient
\[
 \nabla v_t=t^{-1}\mathbf1_{D_t}\nabla U
 \quad\text{almost everywhere}.
\]
Together with \eqref{eq:5.6}, this yields
\begin{equation}
 \int_{\mathbb R^d}|\nabla v_t|^2\,dX
 =t^{-2}\int_{D_t}|\nabla U|^2\,dX\leq\frac Mt.
 \label{eq:5.7}
\end{equation}
Only global square integrability of the gradient is required here; we do not require $v_t\in L^2(\mathbb R^d)$.

We now approximate $v_t$ by functions admissible in \eqref{eq:1.6}. Choose $\chi\in C_c^\infty(\mathbb R^d)$ with $0\leq\chi\leq1$, $\chi=1$ on $B_1$, and $\chi=0$ on $\mathbb R^d\setminus B_2$. Set $\chi_R(X)=\chi(X/R)$ and $w_R=\chi_Rv_t$. Then $w_R\in W^{1,2}(\mathbb R^d)$ has compact support, and
\[
 \nabla w_R-\nabla v_t=(\chi_R-1)\nabla v_t+v_t\nabla\chi_R.
\]
Write
\[
 A_R=\int_{\mathbb R^d}|1-\chi_R|^2|\nabla v_t|^2\,dX,
 \qquad
 Q_R=\int_{\mathbb R^d}v_t^2|\nabla\chi_R|^2\,dX.
\]
By \eqref{eq:5.7},
\[
 0\leq A_R\leq\int_{\mathbb R^d\setminus B_R}|\nabla v_t|^2\,dX
 \longrightarrow0.
\]
For fixed $t$, \eqref{eq:5.1} gives $v_t(X)=O(|X|^{2-d})$. Also, $|\nabla\chi_R|\leq C/R$, with support contained in $B_{2R}\setminus B_R$. Consequently,
\[
 Q_R\leq\frac{C_t}{R^2}\int_R^{2R}r^{2(2-d)}r^{d-1}\,dr
 \leq C_t'R^{2-d}\longrightarrow0.
\]
The cross term is controlled by the Cauchy--Schwarz inequality, explicitly as
\[
 \|\nabla w_R-\nabla v_t\|_{L^2}^2
 \leq A_R+Q_R+2\sqrt{A_RQ_R}\longrightarrow0.
\]
It follows that
\[
 \int_{\mathbb R^d}|\nabla w_R|^2\,dX
 \longrightarrow\int_{\mathbb R^d}|\nabla v_t|^2\,dX.
\]

Write $E_\delta=\{X:\operatorname{dist}(X,E)\leq\delta\}$. By compactness of $E$ and the fact that $U>t$ in a neighborhood of $E$, we can choose $\delta>0$ and $R_0>0$ such that
\[
 E_{3\delta}\subset\{U>t\}\cap B_{R_0}.
\]
Hence $w_R=1$ on $E_{3\delta}$ whenever $R\geq R_0$. Choose a nonnegative standard mollifier $\varphi\in C_c^\infty(B_1)$ with $\int\varphi=1$, and set $\varphi_\epsilon(X)=\epsilon^{-d}\varphi(X/\epsilon)$. For $0<\epsilon<\delta$, the function
\[
 w_{R,\epsilon}=\varphi_\epsilon*w_R
\]
belongs to $C_c^\infty(\mathbb R^d)$. If $\operatorname{dist}(X,E)<\delta$ and $|Y|<\epsilon$, then $X-Y\in E_{2\delta}$, so
\[
 w_{R,\epsilon}(X)
 =\int\varphi_\epsilon(Y)w_R(X-Y)\,dY
 =\int\varphi_\epsilon(Y)\,dY=1.
\]
Thus $w_{R,\epsilon}$ is identically $1$ in an open neighborhood of $E$, and is admissible in \eqref{eq:1.6}. For fixed $R$, continuity of translations in $L^2$ and the inclusion $\nabla w_R\in L^2$ imply
\[
 \|\nabla w_{R,\epsilon}-\nabla w_R\|_{L^2}
 \leq\int\varphi_\epsilon(Y)
 \|\nabla w_R(\cdot-Y)-\nabla w_R\|_{L^2}\,dY
 \longrightarrow0.
\]
First letting $\epsilon\downarrow0$ and then $R\to\infty$, we obtain
\[
 \operatorname{cap}_{D,d}(E)
 \leq\int_{\mathbb R^d}|\nabla v_t|^2\,dX\leq\frac Mt.
\]
By Sard's theorem, the set of critical values has Lebesgue measure zero. We may therefore choose a regular value $t_j$ in each interval $(1-2^{-j},1-2^{-j-1})$. Then $t_j\uparrow1$, and letting $j\to\infty$ proves \eqref{eq:5.2}.
\end{proof}

\subsection{Smooth outer approximation and continuity of capacity from above}

We first describe the construction of smooth compact neighborhoods. Let $K\subset\mathbb R^n$ be a nonempty compact set, let $\varepsilon>0$, and write
\[
 K_\varepsilon=\{x:\operatorname{dist}(x,K)\leq\varepsilon\}.
\]
Take $\delta=\varepsilon/8$ and a nonnegative smooth mollifier $\psi_\delta$ supported in $B_\delta$, with integral $1$. Set
\[
 \chi=\psi_\delta*\mathbf1_{K_{2\delta}}.
\]
Then $\chi\in C_c^\infty(\mathbb R^n)$, $0\leq\chi\leq1$, $\chi=1$ on $K_\delta$, and $\operatorname{supp}\chi\subset K_{3\delta}\subset K_\varepsilon$. Choose a regular value $t\in(0,1)$ and let $L=\{\chi\geq t\}$. The regular value theorem shows that $\partial L=\{\chi=t\}$ is a compact smooth hypersurface and that $L=\overline{\{\chi>t\}}$. Thus $L$ is the closure of a bounded smooth open set and satisfies
\[
 K\subset\operatorname{int}L\subset L\subset K_\varepsilon.
\]
A compact smooth boundary has only finitely many connected components, so both $L$ and its complement have only finitely many connected components.

We next prove continuity from above directly from the energy definition. Let $K_j\subset\mathbb R^n$ be nonempty compact sets with $K_{j+1}\subset K_j$, and set $K=\bigcap_{j=1}^\infty K_j$. These sets are all contained in the compact set $K_1$, so $K\neq\varnothing$. For $0<s<n$, put $k_s(x,y)=|x-y|^{-s}$. In the planar logarithmic case, put $k_0(x,y)=\log(1/|x-y|)$, taking the value $+\infty$ on the diagonal. The following argument applies to both kinds of kernel. Inclusion of the admissible classes of probability measures gives
\[
 E_s(K_j)\uparrow e_s\leq E_s(K).
\]
If $e_s=+\infty$, then $E_s(K)=+\infty$, and continuity of capacity from above follows immediately. Suppose henceforth that $e_s<+\infty$. For each $j$, choose an equilibrium probability measure $\mu_j$ on $K_j$. Weak compactness of probability measures with a common compact support provides a subsequence $\mu_{j_\ell}$ converging weakly to a probability measure $\mu$ on $K_1$. For fixed $m$ and all sufficiently large $\ell$, we have $\mu_{j_\ell}(K_m)=1$. Since $K_m$ is closed, the Portmanteau inequality for weak convergence gives
\[
 \mu(K_m)\geq\limsup_{\ell\to\infty}\mu_{j_\ell}(K_m)=1.
\]
Therefore $\mu(K)=\lim_{m\to\infty}\mu(K_m)=1$, so $\mu\in\mathcal P(K)$.

By the Stone--Weierstrass theorem, finite sums of the form $\sum_{i=1}^N a_i(x)b_i(y)$, with $a_i,b_i\in C(K_1)$, are dense in $C(K_1\times K_1)$. Thus $\mu_{j_\ell}\rightharpoonup\mu$ implies the weak convergence of product measures
\[
 \mu_{j_\ell}\otimes\mu_{j_\ell}\rightharpoonup\mu\otimes\mu.
\]
For every positive integer $N$, the truncated kernel $k_{s,N}=\min\{k_s,N\}$ is continuous and bounded on $K_1\times K_1$. Hence
\[
 \iint k_{s,N}\,d\mu\,d\mu
 =\lim_{\ell\to\infty}\iint k_{s,N}\,d\mu_{j_\ell}\,d\mu_{j_\ell}
 \leq\lim_{\ell\to\infty}E_s(K_{j_\ell})=e_s.
\]
The Riesz kernel is nonnegative, while the logarithmic kernel has a uniform lower bound on the common compact set $K_1\times K_1$. In the latter case, subtract this lower bound before applying the monotone convergence theorem as $N\to\infty$. In either case, we obtain
\[
 E_s(K)\leq I_s(\mu)\leq e_s\leq E_s(K).
\]
Thus $E_s(K_j)\uparrow E_s(K)$. The definition of capacity gives
\[
 C_s(K_j)\downarrow C_s(K).
\]
This conclusion includes the case of zero limiting capacity. In particular, $K_{1/j}\downarrow K$, so $C_s(K_{1/j})\downarrow C_s(K)$. Monotonicity in $\varepsilon$ then gives $C_s(K_\varepsilon)\to C_s(K)$ as $\varepsilon\downarrow0$.

To show that the outer approximation used for the equivalence of Dirichlet capacity in Section~1 does not depend on the energy lemma, we also verify continuity of Dirichlet capacity from above. Let $d\geq3$, let $A_j\subset\mathbb R^d$ be a decreasing sequence of compact sets, and put $A=\bigcap_j A_j$. Monotonicity gives
\[
 \operatorname{cap}_{D,d}(A)
 \leq\lim_{j\to\infty}\operatorname{cap}_{D,d}(A_j).
\]
Given $\epsilon>0$, choose an admissible function $v\in C_c^\infty$ according to \eqref{eq:1.6}, such that
\[
 \int|\nabla v|^2\leq\operatorname{cap}_{D,d}(A)+\epsilon,
 \qquad v\geq1\quad\text{on an open set }V\supset A.
\]
There exists $j_0$ such that $A_j\subset V$ for $j\geq j_0$. Otherwise, choose $x_j\in A_j\setminus V$ and use compactness of $A_1$ to extract a convergent subsequence. Its limit lies in $A\setminus V$, contradicting $A\subset V$. Hence the same $v$ is admissible for $A_j$ for all sufficiently large $j$, and consequently
\[
 \lim_{j\to\infty}\operatorname{cap}_{D,d}(A_j)
 \leq\int|\nabla v|^2
 \leq\operatorname{cap}_{D,d}(A)+\epsilon.
\]
Letting $\epsilon\downarrow0$ proves continuity from above. If the intersection is empty, compactness ensures that some $A_j$ is already empty, and the conclusion remains valid. Together with the preceding construction of smooth neighborhoods, this allows us to choose smooth compact neighborhoods $L_j$ recursively so that
\[
 A\subset\operatorname{int}L_{j+1}\subset L_{j+1}
 \subset\operatorname{int}L_j,\qquad L_j\subset A_{1/j}.
\]
More precisely, at each step construct the next neighborhood inside the open neighborhood $\operatorname{int}L_j\cap\{\operatorname{dist}(\cdot,A)<1/(j+1)\}$ of $A$. Then $L_j\downarrow A$. The equivalence of Dirichlet and Newton capacities on smooth conductors therefore extends, by continuity of both capacities from above, to arbitrary compact sets. In particular, it applies to thin sets embedded in a space of one higher dimension.

Finally, we show that filling holes does not change the capacity with the lower exponent. Let $\Omega_\infty$ be the unbounded connected component of $\mathbb R^n\setminus L$, and set
\[
 F=\mathbb R^n\setminus\Omega_\infty.
\]
This amounts to adjoining to $L$ all bounded components of its complement. Since $L$ is bounded, $F$ is still compact. The boundary $\partial F$ consists of some entire connected components of $\partial L$, so $F$ is the closure of a bounded smooth open set and has connected complement. In particular,
\[
 L\subset F,\qquad\partial F\subset L.
\]
The set $F$ has nonempty interior. The capacity of a ball and monotonicity therefore give $c=C_{n-2}(F)>0$. If $n\geq3$, put $p=n-2$ and take the flux measure $\eta$ constructed for $F$ in Section~2. By \eqref{eq:2.2}, \eqref{eq:2.5}, and \eqref{eq:2.3}, the measure $\nu_F=\eta/c^p$ is a probability measure supported on $\partial F$, and its potential is identically $c^{-p}$ on $F$. Thus
\[
 I_p(\nu_F)=c^{-p}=E_p(F).
\]
If $n=2$, take the probability measure $\nu_F=\mu$ from \eqref{eq:2.9} and \eqref{eq:2.10}. Since $g=0$ on $F$, the logarithmic potential representation \eqref{eq:2.11} similarly gives
\[
 I_0(\nu_F)=-\log c=E_0(F).
\]
In both cases, $\nu_F$ is supported on $\partial F\subset L$, so it is also an admissible probability measure on $L$. Together with monotonicity under $L\subset F$, this yields
\[
 E_{n-2}(F)\leq E_{n-2}(L)\leq I_{n-2}(\nu_F)=E_{n-2}(F).
\]
Consequently,
\[
 C_{n-2}(F)=C_{n-2}(L).
\]
We do not require the filled set $F$ to be contained in $K_\varepsilon$. The approximation below uses only $L\subset K_\varepsilon$ and this equality of capacities.

\subsection{Proof of the main theorem}

\begin{proof}[Proof of Theorem~1.1]
First let $F$ be a smooth set as in Section~2. Apply Lemma~5.1 to $U$ with $d=n+1$. By \eqref{eq:4.8},
\[
 \operatorname{cap}_{D,n+1}(F)\leq(n-1)|S^n|a_nc^{n-1}.
\]
On the other hand, \eqref{eq:1.7} gives
\[
 \operatorname{cap}_{D,n+1}(F)=(n-1)|S^n|C_{n-1}(F)^{n-1}.
\]
The common factor $(n-1)|S^n|$ is strictly positive. Dividing both sides by this factor, we obtain
\begin{equation}
 C_{n-1}(F)^{n-1}\leq a_nC_{n-2}(F)^{n-1}.
 \label{eq:5.8}
\end{equation}
Since $n-1>0$, taking the $(n-1)$st root gives $C_{n-1}(F)\leq a_n^{1/(n-1)}C_{n-2}(F)$, and $a_n^{1/(n-1)}=\kappa_n$.

For an arbitrary nonempty compact set $K$, choose a smooth compact neighborhood $K\subset L\subset K_\varepsilon$, where $K_\varepsilon=\{\operatorname{dist}(x,K)\leq\varepsilon\}$. Fill all bounded components of the complement of $L$, and denote the resulting set by $F$. The equality of capacities under filling proved in the preceding subsection gives $C_{n-2}(F)=C_{n-2}(L)$. By monotonicity and \eqref{eq:5.8},
\begin{equation}
 \begin{aligned}
 C_{n-1}(K)&\leq C_{n-1}(F)
 \leq a_n^{1/(n-1)}C_{n-2}(F)\\
 &=a_n^{1/(n-1)}C_{n-2}(L)
 \leq a_n^{1/(n-1)}C_{n-2}(K_\varepsilon).
 \end{aligned}
 \label{eq:5.9}
\end{equation}
Let $\varepsilon\downarrow0$ and use continuity of capacity from above to obtain the conclusion. Both the smooth approximation and the continuity used here were proved in the preceding subsection. This argument also covers $C_{n-2}(K)=0$; the case of the empty set is immediate.

Finally, by Theorem~1.5, the closed unit ball $B^n$ satisfies
\begin{equation}
 C_{n-2}(B^n)=1,\qquad
 C_{n-1}(B^n)^{n-1}
 =\frac{\Gamma(n/2)}{\sqrt\pi\,\Gamma((n+1)/2)}=a_n.
 \label{eq:5.10}
\end{equation}
Since $a_n=\kappa_n^{n-1}$, this shows that equality holds for the closed unit ball. Translation invariance and $C_s(rK)=rC_s(K)$ then give equality for every closed ball of positive radius. If $\kappa_n$ were decreased, the closed unit ball would no longer satisfy \eqref{eq:1.3}. Thus the constant is optimal.
\end{proof}

\section{Weighted elliptic equations and general exponents}

For a general exponent $q$, we use the same extension function $\rho$, but replace the Laplace operator by the weighted elliptic operator associated with the Riesz kernel. After constructing a supersolution, we compare it pointwise with the normalized $q$-equilibrium potential, and then compare their leading coefficients at infinity.

Fix an integer $n\geq2$ and a real number $q$ satisfying $n-2<q<n$. In particular, $q>0$. Set
\begin{equation}
\beta=q+2-n\in(0,2),\qquad
L_q=\Delta_x+\partial_z^2+\frac{q+1-n}{z}\partial_z
\quad(z\neq0).
\label{eq:6.1}
\end{equation}
Here $X=(x,z)\in\mathbb R^n\times\mathbb R$, and $\Delta_x$ differentiates only with respect to $x$. If $a=q+1-n=\beta-1$, then, for $z\neq0$,
\begin{align*}
|z|^{-a}\operatorname{div}(|z|^a\nabla v)
&=\Delta_xv+\partial_z^2v+|z|^{-a}\frac{d}{dz}|z|^a\,\partial_zv\\
&=\Delta_xv+\partial_z^2v+\frac{a}{z}\partial_zv=L_qv.
\end{align*}
Thus the operator can also be written as
\[
L_q=|z|^{n-q-1}\operatorname{div}(|z|^{q+1-n}\nabla\,\cdot\,).
\]
When $q=n-1$, we have $a=0$, so $L_q=\Delta$.

We next check the integral kernel associated with this operator. Fix $y\in\mathbb R^n$ and put $R_y=|X-(y,0)|$. For $f(r)=r^{-q}$, we have $f'(r)=-qr^{-q-1}$ and $f''(r)=q(q+1)r^{-q-2}$. The formula for the Laplacian of a radial function in $n+1$ dimensions gives
\[
\Delta R_y^{-q}=f''(R_y)+\frac{n}{R_y}f'(R_y)
=q(q+1-n)R_y^{-q-2}.
\]
On the other hand,
\[
\frac{q+1-n}{z}\partial_zR_y^{-q}
=\frac{q+1-n}{z}\bigl(-qzR_y^{-q-2}\bigr)
=-q(q+1-n)R_y^{-q-2}.
\]
Adding these identities yields
\[
L_qR_y^{-q}=0\qquad(z\neq0).
\]

\subsection{Construction of the one-dimensional function and its composition}

First assume that $F$ satisfies the smoothness assumptions of Section~2, let $c=C_{n-2}(F)$, and use the function $\rho$ constructed there. By (2.15) and the identity $1+(q+1-n)=\beta$, we have
\begin{equation}
L_q\rho=\Delta\rho+\frac{q+1-n}{z}\rho_z
=\frac{n-1}{\rho}|\nabla\rho|^2+\beta\frac{\rho_z}{z}
\qquad(z\neq0).
\label{eq:6.2}
\end{equation}
Consequently, for any twice continuously differentiable function $\Phi$ on $(1,\infty)$, the chain rule gives
\[
L_q\Phi(\rho)
=\left(\Phi''(\rho)+\frac{n-1}{\rho}\Phi'(\rho)\right)|\nabla\rho|^2
+\beta\Phi'(\rho)\frac{\rho_z}{z}.
\]
To apply (3.1), consider the boundary value problem
\begin{equation}
\Phi_{n,q}''(r)+\left(\frac{n-1}{r}+\frac{\beta r}{r^2-1}\right)\Phi_{n,q}'(r)=0,
\qquad \Phi_{n,q}(1)=1,\quad \Phi_{n,q}(\infty)=0.
\label{eq:6.3}
\end{equation}
The last condition means that $\lim_{r\to\infty}\Phi_{n,q}(r)=0$. Set
\[
h(t)=t^{1-n}(t^2-1)^{-\beta/2},\qquad t>1.
\]
Logarithmic differentiation gives
\[
\frac{h'(t)}{h(t)}=\frac{1-n}{t}-\frac{\beta t}{t^2-1}.
\]
Near $t=1$, the ratio of $h(t)$ to $(t-1)^{-\beta/2}$ tends to a finite positive number. Since $\beta<2$, the function is integrable at that endpoint. At infinity,
\[
h(t)=t^{1-n-\beta}(1-t^{-2})^{-\beta/2}
=t^{-q-1}(1+O(t^{-2})),
\]
where we used $n+\beta=q+2$. Since $q>0$, the function is also integrable at infinity. Hence the following normalizing constant is finite and strictly positive:
\begin{equation}
Z_{n,q}=\int_1^\infty t^{1-n}(t^2-1)^{-\beta/2}\,dt
=\frac12 B\left(\frac q2,\frac{n-q}{2}\right),
\label{eq:6.4}
\end{equation}
and we define
\begin{equation}
\Phi_{n,q}(r)=\frac1{Z_{n,q}}\int_r^\infty t^{1-n}(t^2-1)^{-\beta/2}\,dt,
\qquad r\geq1.
\label{eq:6.5}
\end{equation}
Here $B$ is the beta function defined in the proof of Lemma~4.1. To check the parameters in \eqref{eq:6.4}, make the substitution $s=t^{-2}$. Then $t=s^{-1/2}$ and $dt=-\frac12s^{-3/2}\,ds$, so
\begin{align*}
Z_{n,q}
&=\frac12\int_0^1s^{(n-1)/2}s^{\beta/2}s^{-3/2}(1-s)^{-\beta/2}\,ds\\
&=\frac12\int_0^1s^{q/2-1}(1-s)^{(n-q)/2-1}\,ds.
\end{align*}
This is the asserted beta integral.

The integral definition shows that $\Phi_{n,q}$ is continuous on $[1,\infty)$ and smooth on $(1,\infty)$, with
\[
\Phi_{n,q}'(r)=-\frac{h(r)}{Z_{n,q}}<0.
\]
Differentiating once more gives $\Phi_{n,q}''(r)=-h'(r)/Z_{n,q}$, and hence
\begin{align*}
&\Phi_{n,q}''(r)+\left(\frac{n-1}{r}+\frac{\beta r}{r^2-1}\right)\Phi_{n,q}'(r)\\
&\qquad=-\frac{h(r)}{Z_{n,q}}
\left(\frac{h'(r)}{h(r)}+\frac{n-1}{r}+\frac{\beta r}{r^2-1}\right)=0.
\end{align*}
The integration limits and integrability give the two boundary conditions in \eqref{eq:6.3}. Substituting this differential equation into the formula for the composition yields
\begin{equation}
L_q\Phi_{n,q}(\rho)
=\beta\Phi_{n,q}'(\rho)
\left(\frac{\rho_z}{z}-\frac{\rho|\nabla\rho|^2}{\rho^2-1}\right)
\leq0\qquad(z\neq0).
\label{eq:6.6}
\end{equation}
Indeed, $\rho>1$ outside $F$, and (3.1) shows that the expression in parentheses is nonnegative. The other two factors satisfy $\beta>0$ and $\Phi_{n,q}'<0$, which gives the stated direction of the inequality.

Let $U_q=\Phi_{n,q}(\rho)$. The continuity and range of $\rho$ imply
\[
U_q=1\quad\text{on }F,\qquad
0<U_q<1\quad\text{on }\mathbb R^{n+1}\setminus F.
\]
The function $U_q$ is continuous on the whole space, smooth outside $F$, and satisfies $U_q(x,z)=U_q(x,-z)$. Integrating the expansion of $h(t)$ at infinity gives
\[
\Phi_{n,q}(r)=\frac1{qZ_{n,q}}r^{-q}+O(r^{-q-2}).
\]
Since $\rho(X)=|X|/c+O(1)$, we also have $\rho(X)^{-q}=c^q|X|^{-q}+O(|X|^{-q-1})$. Therefore
\begin{equation}
U_q(X)=A_{n,q}c^q|X|^{-q}+O(|X|^{-q-1}),
\qquad A_{n,q}=\frac1{qZ_{n,q}}.
\label{eq:6.7}
\end{equation}
The error is uniform in $X/|X|$. Finally, using
$B(u,v)=\Gamma(u)\Gamma(v)/\Gamma(u+v)$ and
$\Gamma(1+q/2)=(q/2)\Gamma(q/2)$, we obtain
\[
A_{n,q}=\frac{2\Gamma(n/2)}{q\Gamma(q/2)\Gamma((n-q)/2)}
=\frac{\Gamma(n/2)}{\Gamma(1+q/2)\Gamma((n-q)/2)}.
\]

\subsection{A comparison principle for functions even across the horizontal hyperplane}

\begin{lemma}
Let $F\subset\mathbb R^n\times\{0\}$ be compact, and let $R>0$ satisfy $F\subset B_R^{n+1}$, where $B_R^{n+1}$ is the open ball of radius $R$ centered at the origin. Put $\Omega=B_R^{n+1}\setminus F$. Suppose that $G$ is bounded on $\Omega$, belongs to $C^2(\Omega)$, and satisfies $G(x,z)=G(x,-z)$. The $C^2$ assumption includes points with $z=0$ and $x\notin F$. If $L_qG\leq0$ on $\Omega\cap\{z\neq0\}$ and
\[
\liminf_{\substack{X\to Y\\X\in\Omega}}G(X)\geq0
\qquad\text{for every }Y\in\partial\Omega,
\]
then $G\geq0$ on $\Omega$.
\end{lemma}

\begin{proof}
We first explain the interpretation of the operator on the hyperplane. Fix $(x_0,0)\in\Omega$. In a small ball about this point, $G$ is twice continuously differentiable and even in $z$, so $G_z(x,0)=0$. The fundamental theorem of calculus gives
\[
\frac{G_z(x,z)}{z}=\int_0^1G_{zz}(x,tz)\,dt\qquad(z\neq0).
\]
Letting $(x,z)\to(x_0,0)$ and using the continuity of $G_{zz}$, we obtain
\[
\frac{G_z(x,z)}{z}\longrightarrow G_{zz}(x_0,0).
\]
Thus the continuous extension of $L_qG$ to this point of the hyperplane is
\begin{align*}
L_qG(x_0,0)
&=\Delta_xG(x_0,0)+(1+q+1-n)G_{zz}(x_0,0)\\
&=\Delta_xG(x_0,0)+\beta G_{zz}(x_0,0).
\end{align*}
Taking a limit in the inequality for $z\neq0$ shows that this value is also nonpositive.

Suppose, to the contrary, that $G(X_0)<0$ for some $X_0\in\Omega$. Since $|X_0|<R$, we may choose $\varepsilon>0$ such that
\[
\varepsilon(R^2-|X_0|^2)<-G(X_0).
\]
Set
\[
G_\varepsilon(X)=G(X)+\varepsilon(R^2-|X|^2).
\]
Then $G_\varepsilon(X_0)<0$. The added term vanishes on the sphere and is nonnegative on $F$, so the lower limit of $G_\varepsilon$ at every boundary point remains nonnegative. Since $G_\varepsilon$ is bounded, $m=\inf_\Omega G_\varepsilon$ is finite and strictly negative. Choose $X_j\in\Omega$ such that $G_\varepsilon(X_j)\to m$. Compactness of the closed ball gives a convergent subsequence, whose limit we denote by $X_*$. If $X_*\in\partial\Omega$, the boundary lower-limit condition would imply $m\geq0$, a contradiction. Hence $X_*\in\Omega$, and continuity gives $G_\varepsilon(X_*)=m$.

On the other hand, direct computation yields
\[
L_q|X|^2=2n+2+2(q+1-n)=2(q+2),
\]
and the continuous interpretation on the hyperplane gives the same value. Thus
\[
L_qG_\varepsilon\leq-2\varepsilon(q+2)<0
\qquad\text{throughout }\Omega.
\]
At the interior minimum $X_*$, we have $\nabla G_\varepsilon=0$, and the second derivative in each coordinate direction is nonnegative. If $z_*\neq0$, the drift term vanishes, so
$L_qG_\varepsilon(X_*)=\Delta G_\varepsilon(X_*)\geq0$. If $z_*=0$, then $\beta>0$ gives
\[
L_qG_\varepsilon(X_*)=\Delta_xG_\varepsilon(X_*)
+\beta(G_\varepsilon)_{zz}(X_*)\geq0.
\]
Both cases contradict the strict negative inequality. Therefore $G\geq0$.
\end{proof}

\begin{proof}[Proof of Theorem~1.2]
First suppose that $F$ is the closure of a nonempty bounded smooth open set and has connected complement. We begin by checking that $C_q(F)>0$. Let $m$ denote the $n$-dimensional Lebesgue measure of $F$, and put
$D=\operatorname{diam}(F)=\sup\{|x-y|:x,y\in F\}$.
Since $F$ has nonempty interior, $m>0$ and $D>0$. Let $\tau$ be normalized volume measure on $F$, that is,
$d\tau(x)=m^{-1}\mathbf1_F(x)\,dx$. Since $q<n$,
\begin{align*}
I_q(\tau)
&=\frac1{m^2}\int_F\int_F|x-y|^{-q}\,dy\,dx\\
&\leq\frac1m\int_{B_D^n}|h|^{-q}\,dh
=\frac{|S^{n-1}|D^{n-q}}{(n-q)m}<\infty.
\end{align*}
Hence $E_q(F)<\infty$. Moreover, $|x-y|\leq D$ for $x,y\in F$, so the energy of every probability measure is at least $D^{-q}$, and therefore $E_q(F)>0$. In particular,
$C_q(F)^q=E_q(F)^{-1}$ is finite and strictly positive. The distances in these integrals are the same as after embedding $F$ in $\mathbb R^{n+1}$; thus the fact that $F$ has zero volume in the higher-dimensional space does not affect this finite-energy conclusion.

By Theorem~1.3, there exists a $q$-equilibrium probability measure $\sigma_q$ on $F$. Define the normalized $q$-equilibrium potential by
\[
W_q(X)=C_q(F)^q\int_F|X-(y,0)|^{-q}\,d\sigma_q(y).
\]
The upper bound for the equilibrium potential used here corresponds to the parameter $\alpha=n-q\in(0,2)$. Thus Theorem~1.3 gives
\[
\int_F|x-y|^{-q}\,d\sigma_q(y)\leq E_q(F)
\qquad\text{for every }x\in\mathbb R^n.
\]
Since $q>0$ and $|X-(y,0)|\geq|x-y|$, we also have
\begin{align*}
0\leq W_q(x,z)
&\leq C_q(F)^q\int_F|x-y|^{-q}\,d\sigma_q(y)\\
&\leq C_q(F)^qE_q(F)=1.
\end{align*}
We use the upper bound on the whole space; we do not require the equilibrium potential to equal the equilibrium energy at every point of $F$.

In any compact neighborhood disjoint from $F$, the distance between $X$ and the integration point $(y,0)$ has a strictly positive lower bound. Hence the kernel and its derivatives of any fixed order are uniformly bounded in $y\in F$ on that neighborhood, and repeated differentiation under the integral sign is justified. This proves that $W_q$ is smooth outside $F$ and satisfies $W_q(x,z)=W_q(x,-z)$. By the kernel calculation at the beginning of this section,
\[
L_qW_q=0\qquad(z\neq0).
\]
Since $F$ is bounded, set $M=\sup_{y\in F}|y|$. For $r=|X|$,
\[
|X-(y,0)|^{-q}
=r^{-q}\left(1-\frac{2X\cdot(y,0)}{r^2}+\frac{|y|^2}{r^2}\right)^{-q/2}
=r^{-q}(1+O(r^{-1})).
\]
The last error is uniform in $y\in F$ and $X/r$. Indeed, the absolute value of the terms other than $1$ in parentheses is at most $2M/r+M^2/r^2$, and the function $(1+s)^{-q/2}$ has bounded derivative near $s=0$. Using $\sigma_q(F)=1$, we obtain
\[
W_q(X)=C_q(F)^q r^{-q}+O(r^{-q-1}).
\]

Let $H=U_q-W_q$. The ranges established above give $-1\leq H\leq1$. The function is smooth outside $F$, even in $z$, and satisfies
\[
L_qH\leq0\qquad(z\neq0).
\]
For every $Y\in F$, continuity of $U_q$ gives $U_q(X)\to1$. Since $W_q\leq1$, it follows that
\[
\liminf_{\substack{X\to Y\\X\notin F}}H(X)
\geq\liminf_{\substack{X\to Y\\X\notin F}}(U_q(X)-1)=0.
\]
The far-field expansions of the two functions also imply $H(X)\to0$, uniformly in direction.

Fix $X_0\notin F$ and $\delta>0$. By the limit at infinity, choose $R_0$ so that $|H(X)|\leq\delta$ for all $|X|\geq R_0$. Then choose
\[
R>\max\left\{R_0,|X_0|,\sup_{Y\in F}|Y|\right\}.
\]
We have $H+\delta\geq0$ on $\partial B_R^{n+1}$ and
$\liminf(H+\delta)\geq\delta>0$ on $F$. Applying Lemma~6.1 to $H+\delta$ in $B_R^{n+1}\setminus F$ gives $H(X_0)\geq-\delta$. With $X_0$ fixed and $\delta>0$ arbitrary, we let $\delta\downarrow0$ to obtain $H(X_0)\geq0$. Since $X_0\notin F$ was arbitrary,
\[
U_q(X)\geq W_q(X)\qquad(X\notin F).
\]

We now compare the far-field coefficients. Take $X=(0,r)\in\mathbb R^n\times\mathbb R$, with $r>0$. Multiplying the nonnegative difference above by $r^q$ and using the two far-field expansions gives
\[
0\leq r^q\bigl(U_q(0,r)-W_q(0,r)\bigr)
=A_{n,q}c^q-C_q(F)^q+O(r^{-1}).
\]
Letting $r\to\infty$, we obtain
\[
C_q(F)^q\leq A_{n,q}c^q=A_{n,q}C_{n-2}(F)^q.
\]
Since $q>0$, taking the $q$th root yields
$C_q(F)\leq A_{n,q}^{1/q}C_{n-2}(F)$.

Finally, let $K\subset\mathbb R^n$ be an arbitrary nonempty compact set. For $\varepsilon>0$, set
$K_\varepsilon=\{x:\operatorname{dist}(x,K)\leq\varepsilon\}$.
By the smooth outer approximation of Section~5.2, there exists a smooth compact neighborhood $L$ such that $K\subset L\subset K_\varepsilon$. Fill in the bounded components of the complement of $L$ to obtain a smooth set $F$ with connected complement and $L\subset F$. Section~5.2 proves that $C_{n-2}(F)=C_{n-2}(L)$. Monotonicity of capacity and the smooth case just proved give, successively,
\begin{align*}
C_q(K)&\leq C_q(L)\leq C_q(F)\\
&\leq A_{n,q}^{1/q}C_{n-2}(F)
=A_{n,q}^{1/q}C_{n-2}(L)\\
&\leq A_{n,q}^{1/q}C_{n-2}(K_\varepsilon).
\end{align*}
Take $\varepsilon=1/j$. Then $K_{1/j}$ decreases to $K$. By the continuity under decreasing intersections proved in Section~5.2,
$C_{n-2}(K_{1/j})\downarrow C_{n-2}(K)$. Letting $j\to\infty$, we conclude that
\[
C_q(K)\leq A_{n,q}^{1/q}C_{n-2}(K).
\]
This limit also covers the case $C_{n-2}(K)=0$. Both capacities of the empty set are zero by definition, so the conclusion also holds for the empty set.

By Theorem~1.5, the closed unit ball $B^n$ satisfies
\[
C_{n-2}(B^n)=1,\qquad
C_q(B^n)^q=\frac{\Gamma(n/2)}{\Gamma(1+q/2)\Gamma((n-q)/2)}=A_{n,q}.
\]
Thus equality holds for the closed unit ball. Capacity is translation invariant and satisfies $C_s(rK)=rC_s(K)$, including when $s=0$. Hence every closed ball of positive radius attains equality, and the constant $A_{n,q}^{1/q}$ cannot be decreased.
\end{proof}


\clearpage
\begin{thebibliography}{8}
\bibitem{ref1}
G. P\'olya and G. Szeg\H{o},
\emph{Inequalities for the capacity of a condenser},
American Journal of Mathematics \textbf{67} (1945), no.~1, 1--32.
\href{https://doi.org/10.2307/2371912}{doi:10.2307/2371912}.

\bibitem{ref2}
C. Clark and R. S. Laugesen,
\emph{Toward P\'olya and Szeg\H{o}'s conjecture for logarithmic vs Newtonian capacity},
preprint, September 28, 2026.
\href{https://arxiv.org/abs/2609.35438v1}{arXiv:2609.35438v1}.

\bibitem{ref3}
C. Clark and R. S. Laugesen,
\emph{Maximizing Riesz capacity ratios: conjectures and theorems},
SIAM Journal on Mathematical Analysis \textbf{57} (2025), no.~5, 5683--5717.
\href{https://doi.org/10.1137/24M171992X}{doi:10.1137/24M171992X}.

\bibitem{ref4}
C. Clark and R. S. Laugesen,
\emph{Riesz capacity: monotonicity, continuity, diameter and volume},
Analysis and Mathematical Physics \textbf{15} (2025), article~4.
\href{https://doi.org/10.1007/s13324-024-01000-2}{doi:10.1007/s13324-024-01000-2}.

\bibitem{ref5}
N. S. Landkof,
\emph{Foundations of Modern Potential Theory},
Grundlehren der mathematischen Wissenschaften, vol.~180, Springer, 1972.
\href{https://doi.org/10.1007/978-3-642-65183-0}{doi:10.1007/978-3-642-65183-0}.

\bibitem{ref6}
S. V. Borodachov, D. P. Hardin and E. B. Saff,
\emph{Discrete Energy on Rectifiable Sets},
Springer Monographs in Mathematics, Springer, 2019.
\href{https://doi.org/10.1007/978-0-387-84808-2}{doi:10.1007/978-0-387-84808-2}.

\bibitem{ref7}
E. B. Saff,
\emph{Logarithmic potential theory with applications to approximation theory},
Surveys in Approximation Theory \textbf{5} (2010), 165--200.
\href{https://arxiv.org/abs/1010.3760}{arXiv:1010.3760}.

\bibitem{ref8}
I. E. Pritsker, E. B. Saff and W. Wise,
\emph{Reverse triangle inequalities for Riesz potentials and connections with polarization},
Journal of Mathematical Analysis and Applications \textbf{410} (2014), no.~2, 868--881.
\href{https://doi.org/10.1016/j.jmaa.2013.09.002}{doi:10.1016/j.jmaa.2013.09.002}.
\end{thebibliography}
\end{document}